\documentclass[11pt]{article}
\usepackage[margin=1in]{geometry}
\usepackage[T1]{fontenc}
\usepackage{lmodern}
\usepackage{microtype}
\usepackage{amsmath,amssymb,amsthm,mathtools}
\usepackage{algorithm}
\usepackage{algpseudocode}
\usepackage{booktabs,array}
\usepackage[authoryear,square]{natbib}
\usepackage{xcolor}
\usepackage{xurl}
\usepackage[colorlinks=true,linkcolor=blue!55!black,citecolor=blue!55!black,urlcolor=blue!55!black]{hyperref}
\usepackage{enumitem}
\usepackage{float}
\usepackage{needspace}
\hypersetup{pdftitle={Variance Reduction with Delayed Feedback for Non-Convex Optimization}}
\setlist[itemize]{topsep=4pt,itemsep=2pt}
\allowdisplaybreaks[2]
\newtheorem{theorem}{Theorem}
\newtheorem{lemma}{Lemma}
\newtheorem{corollary}{Corollary}

\newtheorem{assumption}{Assumption}
\theoremstyle{remark}
\newtheorem{remark}{Remark}
\newcommand{\E}{\mathbb E}
\newcommand{\R}{\mathbb R}
\newcommand{\norm}[1]{\left\|#1\right\|}
\newcommand{\inner}[2]{\left\langle#1,#2\right\rangle}
\newcommand{\ind}{\mathbf 1}

\title{Variance Reduction with Delayed Feedback\newline for Non-Convex Optimization}
\author{
{Wei Jiang}$^{1,2}$ \quad
{Dingzhi Yu}$^{2}$ \quad
{Zechao Li}$^{1}$ \\[0.5em]
\small$^1$School of Computer Science and Engineering, 
Nanjing University of Science and Technology\\
\small$^2$National Key Laboratory for Novel Software Technology, Nanjing University\\
\small jiangw@lamda.nju.edu.cn,\quad yudz@lamda.nju.edu.cn, \quad zechao.li@njust.edu.cn
}
\begin{document}
\maketitle
\vspace{-1.5em}
\begin{abstract}
We study variance reduction for non-convex optimization when gradient evaluations return with a delay. The server can update the parameter while earlier evaluations remain unfinished. We first analyze an ordered recursive momentum estimator for stochastic optimization. Under average smoothness and bounded variance, its oracle complexity is $O(1+\sigma^2\epsilon^{-2}+L\Delta_f(1+\sigma/\epsilon+\tau)\epsilon^{-2})$, where the target is an expected squared gradient norm of at most $\epsilon^2$ and $\tau$ bounds the feedback delay. We next study finite sums with $n$ components and obtain $O(n+L\Delta_f(\sqrt n+\min\{\tau,n-1\})\epsilon^{-2})$ evaluations. An asynchronous snapshot variant continues updating during full-gradient computations under an explicit bound on their completion time. For the stochastic method, we further replace a maximum-delay bound by a bound on the weighted amount of missing feedback. 
\end{abstract}

\section{Introduction}

This paper studies variance reduction for stochastic non-convex optimization with delayed feedback. Recursive estimators such as SARAH, SPIDER, and STORM employ the same sample at two nearby points to reduce the noise in a gradient difference \citep{nguyen2017sarah,fang2018spider,cutkosky2019storm,wang2019spiderboost}. With delayed evaluations, however, the correction associated with one update may be unavailable when the next update is performed. Several corrections can be outstanding, and their arrival order can depend on the sampled data.

A direct approach fixes the query points and waits for a complete batch before updating the parameter. Waiting preserves the usual variance-reduction oracle complexity, since it does not consume gradient evaluations. We investigate methods that continue to update the parameter while corrections are outstanding. Each correction is evaluated at two consecutive server iterates, and the server incorporates it using its generation index when it returns. We count all issued evaluations, including unfinished tasks.

We first consider stochastic optimization under bounded variance and average smoothness. After one initialization batch, the method issues one two-point task per update without further batch-completion waits. To obtain $\E\norm{\nabla f(x_R)}^2\leq\epsilon^2$, its oracle complexity is
\[
O\!\left(1+\frac{\sigma^2}{\epsilon^2}
+\frac{L\Delta_f}{\epsilon^2}
\left(1+\frac{\sigma}{\epsilon}+\tau\right)\right).
\]
Here $\tau$ bounds the correction delay in server updates, $\sigma^2$ bounds the oracle variance, $L$ is the average-smoothness constant, and $\Delta_f$ is the initial objective gap. We next consider a finite sum of $n$ components. Using exact gradient refreshes every $n$ updates and one sampled difference per update, we obtain
\[
O\!\left(n+\frac{L\Delta_f(\sqrt n+\min\{\tau,n-1\})}{\epsilon^2}\right).
\]
Thus the delay enters additively in both bounds. In the finite-sum result, a refresh also limits how long an unfinished correction can affect the estimator.

We further study two extensions of these results. For stochastic optimization, we control the outstanding corrections through their exponentially decaying weights in Section~\ref{sec:outstanding}. The resulting condition allows arbitrarily late responses and repeated permanent losses, provided that these losses are sufficiently separated relative to the momentum parameter. For finite sums, Section~\ref{sec:snapshots} removes the refresh waits after initialization by computing exact gradients at fixed snapshot points while the server continues to update. Under an additional bound $D$ on the completion delay of each exact snapshot, this method has complexity
\[
O\!\left(n+\frac{L\Delta_f(\sqrt n+\sqrt D+\tau)}{\epsilon^2}\right).
\]
The service required for these exact snapshots is stated separately from the delay assumption for sampled corrections.

\section{Problem Setting}
\label{sec:setting}

We consider the non-convex optimization problem
\begin{equation}
\min_{x\in\R^d} f(x).
\label{eq:problem}
\end{equation}
The objective is in the form of $f(x)=\E_\xi[f(x;\xi)]$~(expectation) or $f(x)=n^{-1}\sum_{i=1}^n f_i(x)$~(finite-sum). 

\begin{assumption}[Lower boundedness]
\label{ass:lower}
The objective satisfies $f_*:=\inf_x f(x)>-\infty$. The initial point $x_0$ satisfies $f(x_0)-f_*\leq\Delta_f$.
\end{assumption}

Our goal is to return a point $x_R$ such that
\begin{equation}
\E\bigl[\norm{\nabla f(x_R)}^2\bigr]\leq\epsilon^2.
\label{eq:target}
\end{equation}
The output index $R$ is chosen uniformly from the index set $\{0,\ldots,T-1\}$. In particular, equation \eqref{eq:target} also implies $\E[\norm{\nabla f(x_R)}]\leq\epsilon$. All guarantees below are in expectation.

\subsection{Delayed two-point queries}

At update $t$, the server has an estimator $v_t$ and computes
\begin{equation}
x_{t+1}=x_t-\eta v_t.
\label{eq:update}
\end{equation}
It then requests a correction evaluated at the fixed pair $(x_t,x_{t+1})$. A fresh sample is drawn after this pair has been chosen. The two gradient evaluations in the task use the same sample, and we treat such a task as a single response.
Let $A_t$ denote the first server round in which the response to task $t$ is available. Our initial assumption is
\begin{equation}
t+1\leq A_t\leq t+1+\tau,
\qquad t\geq0,
\label{eq:delay}
\end{equation}
where $\tau$ is a nonnegative integer. Thus $\tau=0$ means that the correction is available for the next update. Tasks can arrive out of order. 

At the beginning of a round, the server incorporates all responses available for that round. It performs \eqref{eq:update} even if no new response has arrived, and then issues the next task. We do not assume that each received message triggers exactly one update. Initialization uses a complete preselected batch at a fixed point. The basic finite-sum method also waits for a complete full gradient at each refresh. Section~\ref{sec:snapshots} gives a variant that computes these later full gradients while the server continues updating. Each snapshot still evaluates all $n$ components at its fixed point.

\subsection{The descent inequality}

\begin{lemma}
\label{lem:descent}
Suppose that $f$ is $L$-smooth and \eqref{eq:update} is used with $\eta>0$. Then
\begin{align}
f(x_{t+1})
\leq{}&f(x_t)-\frac{\eta}{2}\norm{\nabla f(x_t)}^2
-\frac{\eta}{2}(1-L\eta)\norm{v_t}^2
+\frac{\eta}{2}\norm{v_t-\nabla f(x_t)}^2.
\label{eq:descent}
\end{align}
Consequently, under Assumption~\ref{ass:lower},
\begin{align}
\sum_{t=0}^{T-1}\E\norm{\nabla f(x_t)}^2
\leq{}&\frac{2\Delta_f}{\eta}
-(1-L\eta)\sum_{t=0}^{T-1}\E\norm{v_t}^2
+\sum_{t=0}^{T-1}\E\norm{v_t-\nabla f(x_t)}^2.
\label{eq:descent-sum}
\end{align}
\end{lemma}

\begin{proof}
Smoothness and \eqref{eq:update} give
\begin{align*}
f(x_{t+1})
&\leq f(x_t)+\inner{\nabla f(x_t)}{x_{t+1}-x_t}
+\frac L2\norm{x_{t+1}-x_t}^2\\
&=f(x_t)-\eta\inner{\nabla f(x_t)}{v_t}
+\frac{L\eta^2}{2}\norm{v_t}^2.
\end{align*}
Since
\[
\norm{v_t-\nabla f(x_t)}^2
=\norm{v_t}^2+\norm{\nabla f(x_t)}^2
-2\inner{\nabla f(x_t)}{v_t},
\]
we obtain
\begin{align*}
f(x_{t+1})
\leq{}&f(x_t)-\frac\eta2\norm{\nabla f(x_t)}^2
-\frac\eta2\norm{v_t}^2
+\frac\eta2\norm{v_t-\nabla f(x_t)}^2
+\frac{L\eta^2}{2}\norm{v_t}^2\\
={}&f(x_t)-\frac\eta2\norm{\nabla f(x_t)}^2
-\frac\eta2(1-L\eta)\norm{v_t}^2
+\frac\eta2\norm{v_t-\nabla f(x_t)}^2,
\end{align*}
which proves \eqref{eq:descent}. Rearranging and taking expectations gives
\begin{align*}
\E\norm{\nabla f(x_t)}^2
\leq{}&\frac2\eta\bigl(\E f(x_t)-\E f(x_{t+1})\bigr)
-(1-L\eta)\E\norm{v_t}^2
+\E\norm{v_t-\nabla f(x_t)}^2.
\end{align*}
The objective differences telescope:
\begin{align*}
\sum_{t=0}^{T-1}\bigl(\E f(x_t)-\E f(x_{t+1})\bigr)
=f(x_0)-\E f(x_T)\leq f(x_0)-f_*\leq\Delta_f.
\end{align*}
Summing the preceding inequalities proves \eqref{eq:descent-sum}. 
\end{proof}

\section{Variance Reduction for Stochastic Optimization}
\label{sec:stochastic}

We first consider $f(x)=\E_\xi[f(x;\xi)]$. An exact full gradient is generally unavailable. We use a STORM estimator, whose momentum term gradually reduces the effect of old estimation errors.

\begin{assumption}[Stochastic oracle]
\label{ass:stochastic}
For all $x,y\in\R^d$, the stochastic gradient satisfies
\begin{align}
\E_\xi[\nabla f(x;\xi)]&=\nabla f(x),\label{eq:unbiased}\\
\E_\xi\norm{\nabla f(x;\xi)-\nabla f(x)}^2&\leq\sigma^2,\label{eq:variance}\\
\E_\xi\norm{\nabla f(x;\xi)-\nabla f(y;\xi)}^2
&\leq L^2\norm{x-y}^2.\label{eq:average-smooth}
\end{align}
Here $\sigma\geq0$, $L>0$, and the same sample $\xi$ is used at both points in \eqref{eq:average-smooth}.
\end{assumption}

These are the usual bounded-variance and average-smoothness conditions for recursive variance reduction. In particular, we do not assume a uniform bound on $\norm{\nabla f(x;\xi)}$. Jensen's inequality applied to \eqref{eq:average-smooth} implies that $f$ is $L$-smooth.

\subsection{The proposed method}

We first draw $B$ fresh samples and wait for all of their evaluations at $x_0$ to compute
\begin{equation}
v_0=\frac{1}{B}\sum_{k=1}^B\nabla f(x_0;\xi_{0,k}).
\label{eq:initialization}
\end{equation}
After the update from $x_t$ to $x_{t+1}$, task $t$ draws a fresh sample $\xi_t$ and returns
\begin{align}
\delta_t
&=\nabla f(x_{t+1};\xi_t)-(1-\alpha)\nabla f(x_t;\xi_t)\nonumber\\
&=\nabla f(x_{t+1};\xi_t)-\nabla f(x_t;\xi_t)
+\alpha\nabla f(x_t;\xi_t),
\label{eq:online-correction}
\end{align}
where $0<\alpha<1$. 
The server uses the available corrections with weights determined by their generation times:
\begin{equation}
v_t=(1-\alpha)^t v_0+
\sum_{j=0}^{t-1}\ind_{\{A_j\leq t\}}(1-\alpha)^{t-j-1}\delta_j.
\label{eq:online-estimator}
\end{equation}
Thus a correction does not receive the weight of a newly generated sample merely because it arrives late. It receives the weight it would have had if it had been available from the next round. Formula \eqref{eq:online-estimator} can be implemented by
\begin{equation}
v_t=(1-\alpha)v_{t-1}
+\sum_{j:A_j=t}(1-\alpha)^{t-j-1}\delta_j,
\qquad t\geq1.
\label{eq:online-recursion}
\end{equation}
If no response arrives in a round, the estimator is multiplied by $1-\alpha$, the parameter is updated, and another task is issued. Multiple responses can be incorporated in one round.
\begin{algorithm}[tb]
\caption{Delayed recursive momentum for stochastic optimization}
\label{alg:online}
\begin{algorithmic}[1]
\Require Initial point $x_0$, update number $T$, initialization size $B$, parameters $\eta,\alpha$
\State Compute $v_0$ using the complete preselected batch in \eqref{eq:initialization}.
\For{$t=0,\ldots,T-1$}
  \If{$t\geq1$}
    \State Update $v_t$ by \eqref{eq:online-recursion}, including all responses arriving in round $t$.
  \EndIf
  \State Update $x_{t+1}=x_t-\eta v_t$.
  \State Draw a fresh $\xi_t$ and issue the two-point task \eqref{eq:online-correction}.
\EndFor
\State Choose $R$ uniformly from $\{0,\ldots,T-1\}$ and return $x_R$.
\end{algorithmic}
\end{algorithm}

\subsection{Convergence guarantee}

\begin{theorem}
\label{thm:online}
Under Assumptions~\ref{ass:lower} and~\ref{ass:stochastic}, suppose that \eqref{eq:delay} holds. For any $\epsilon>0$, choose
\begin{align}
\eta&=\frac{1}{16L(1+\sigma/\epsilon+\tau)},
&\alpha&=\frac{1}{16(1+\sigma/\epsilon)(1+\sigma/\epsilon+\tau)},
\label{eq:online-params-a}\\
T&=\max\left\{1,\left\lceil
\frac{64L\Delta_f(1+\sigma/\epsilon+\tau)}{\epsilon^2}
\right\rceil\right\}.
\label{eq:online-params-b}
\end{align}
Set the initialization size to
\begin{equation}
B=\max\left\{1,\left\lceil
\frac{16\sigma^2}{\epsilon^2}
\min\left\{1,\frac1{\alpha T}\right\}
\right\rceil\right\}.
\label{eq:online-batch}
\end{equation}
Then Algorithm~\ref{alg:online} guarantees
\begin{equation}
\E\norm{\nabla f(x_R)}^2\leq\epsilon^2.
\label{eq:online-target}
\end{equation}
In particular, its oracle complexity is
\begin{equation}
O\!\left(1+\frac{\sigma^2}{\epsilon^2}
+\frac{L\Delta_f}{\epsilon^2}
\left(1+\frac{\sigma}{\epsilon}+\tau\right)\right).
\label{eq:online-complexity}
\end{equation}
\end{theorem}

\begin{remark}
For $0<\epsilon\leq\sigma$, \eqref{eq:online-complexity} becomes
\[
O\!\left(\frac{\sigma^2}{\epsilon^2}
+\frac{L\Delta_f\sigma}{\epsilon^3}
+\frac{L\Delta_f\tau}{\epsilon^2}\right).
\]
Hence $\tau=O(\sigma/\epsilon)$ preserves the usual stochastic variance-reduction order.
\end{remark}

\begin{remark}[Initialization]
The factor $\min\{1,1/(\alpha T)\}$ in \eqref{eq:online-batch} uses the decay of the initial estimation error. For $\Delta_f>0$, the choices above imply
\[
\alpha T\geq\frac{4L\Delta_f}{(1+\sigma/\epsilon)\epsilon^2},
\qquad
B\leq1+\min\left\{
\frac{16\sigma^2}{\epsilon^2},
\frac{4\sigma^2(1+\sigma/\epsilon)}{L\Delta_f}
\right\}.
\]
Thus a longer run can start with a smaller complete batch. This reduces the number of evaluations required before the first update while preserving the complexity in \eqref{eq:online-complexity}.
\end{remark}

\subsection{Estimation error}

We introduce an auxiliary estimator that contains all generated corrections:
\begin{equation}
u_0=v_0,\qquad u_{t+1}=(1-\alpha)u_t+\delta_t.
\label{eq:virtual-estimator}
\end{equation}
The sequence $u_t$ is used only for analysis. Its corrections are evaluated at the actual iterates of Algorithm~\ref{alg:online} and it does not generate a second sequence of parameters. Define
\begin{equation}
\mathcal M_t=\{j:0\leq j<t,\ A_j>t\}.
\label{eq:online-pending}
\end{equation}
For $t\geq2$, repeated substitution in \eqref{eq:virtual-estimator} gives
\begin{align*}
u_t
&=(1-\alpha)u_{t-1}+\delta_{t-1}\\
&=(1-\alpha)^2u_{t-2}+(1-\alpha)\delta_{t-2}+\delta_{t-1}\\
&=\cdots=(1-\alpha)^t v_0
+\sum_{j=0}^{t-1}(1-\alpha)^{t-j-1}\delta_j.
\end{align*}
The last equality also holds for $t=0,1$, with an empty sum when $t=0$. Comparing with \eqref{eq:online-estimator},
\begin{align*}
u_t-v_t
&=\sum_{j=0}^{t-1}\bigl(1-\ind_{\{A_j\leq t\}}\bigr)
(1-\alpha)^{t-j-1}\delta_j\\
&=\sum_{j\in\mathcal M_t}(1-\alpha)^{t-j-1}\delta_j.
\end{align*}
Consequently,
\begin{equation}
v_t-\nabla f(x_t)
=u_t-\nabla f(x_t)-
\sum_{j\in\mathcal M_t}(1-\alpha)^{t-j-1}\delta_j.
\label{eq:online-decomposition}
\end{equation}
The first term includes all generated samples, whereas the second term records missing responses. We bound these two terms separately.

\begin{lemma}
\label{lem:online-virtual}
Under Assumption~\ref{ass:stochastic}, for any $0<\alpha<1$ and $T\geq1$,
\begin{align}
\sum_{t=0}^{T-1}\E\norm{u_t-\nabla f(x_t)}^2
\leq{}&\frac{\sigma^2}{B}\min\left\{T,\frac1\alpha\right\}
+\frac{2L^2}{\alpha}\sum_{t=0}^{T-1}\E\norm{x_{t+1}-x_t}^2
+2\alpha\sigma^2T,
\label{eq:virtual-error}
\end{align}
and each correction satisfies
\begin{equation}
\E\norm{\delta_t}^2
\leq2L^2\E\norm{x_{t+1}-x_t}^2
+2\alpha^2\E\norm{\nabla f(x_t)}^2+2\alpha^2\sigma^2.
\label{eq:online-packet}
\end{equation}
These bounds do not require a maximum-delay assumption.
\end{lemma}

\begin{proof}
Center the new correction at its population value:
\begin{align}
\zeta_t={}&\delta_t-
\bigl(\nabla f(x_{t+1})-(1-\alpha)\nabla f(x_t)\bigr)\nonumber\\
={}&\nabla f(x_{t+1};\xi_t)-\nabla f(x_t;\xi_t)
-\bigl(\nabla f(x_{t+1})-\nabla f(x_t)\bigr)\nonumber\\
&+\alpha\bigl(\nabla f(x_t;\xi_t)-\nabla f(x_t)\bigr).
\label{eq:online-noise-definition}
\end{align}
The points $x_t,x_{t+1}$ are chosen before $\xi_t$. Unbiasedness therefore gives
\begin{align*}
\E_{\xi_t}\zeta_t
={}&\nabla f(x_{t+1})-\nabla f(x_t)
-\bigl(\nabla f(x_{t+1})-\nabla f(x_t)\bigr)\\
&+\alpha\bigl(\nabla f(x_t)-\nabla f(x_t)\bigr)=0,
\end{align*}
where the chosen pair is held fixed in this sample average. Also,
\begin{align*}
&\E_{\xi_t}\bigl\|
\nabla f(x_{t+1};\xi_t)-\nabla f(x_t;\xi_t)
-\bigl(\nabla f(x_{t+1})-\nabla f(x_t)\bigr)\bigr\|^2\\
&\quad=\E_{\xi_t}\norm{\nabla f(x_{t+1};\xi_t)-\nabla f(x_t;\xi_t)}^2
-\norm{\nabla f(x_{t+1})-\nabla f(x_t)}^2\\
&\quad\leq\E_{\xi_t}\norm{\nabla f(x_{t+1};\xi_t)-\nabla f(x_t;\xi_t)}^2.
\end{align*}
Using $\norm{a+b}^2\leq2\norm{a}^2+2\norm{b}^2$ in \eqref{eq:online-noise-definition}, then averaging the preceding randomness, yields
\begin{align}
\E\norm{\zeta_t}^2
&\leq2\E\norm{\nabla f(x_{t+1};\xi_t)-\nabla f(x_t;\xi_t)}^2
+2\alpha^2\sigma^2\nonumber\\
&\leq2L^2\E\norm{x_{t+1}-x_t}^2+2\alpha^2\sigma^2.
\label{eq:online-noise}
\end{align}
Subtracting $\nabla f(x_{t+1})$ from \eqref{eq:virtual-estimator},
\begin{align*}
u_{t+1}-\nabla f(x_{t+1})
&=(1-\alpha)u_t+\delta_t-\nabla f(x_{t+1})\\
&=(1-\alpha)\bigl(u_t-\nabla f(x_t)\bigr)+\zeta_t.
\end{align*}
Thus, for $t\geq2$,
\begin{align*}
u_t-\nabla f(x_t)
&=(1-\alpha)\bigl(u_{t-1}-\nabla f(x_{t-1})\bigr)+\zeta_{t-1}\\
&=(1-\alpha)^2\bigl(u_{t-2}-\nabla f(x_{t-2})\bigr)
+(1-\alpha)\zeta_{t-2}+\zeta_{t-1}.
\end{align*}
Continuing to $u_0=v_0$ gives, for every $t\geq0$,
\begin{equation}
u_t-\nabla f(x_t)
=(1-\alpha)^t\bigl(v_0-\nabla f(x_0)\bigr)
+\sum_{j=0}^{t-1}(1-\alpha)^{t-j-1}\zeta_j.
\label{eq:virtual-unroll}
\end{equation}
For each $k\geq0$, average $\xi_k$ with the preceding samples and its point pair fixed. The zero mean just proved gives
\[
\E\inner{v_0-\nabla f(x_0)}{\zeta_k}=0.
\]
For $0\leq j<k$, the same averaging gives $\E\inner{\zeta_j}{\zeta_k}=0$.
Expanding the square in \eqref{eq:virtual-unroll} now gives
\begin{align*}
\E\norm{u_t-\nabla f(x_t)}^2
={}&(1-\alpha)^{2t}\E\norm{v_0-\nabla f(x_0)}^2
+\sum_{j=0}^{t-1}(1-\alpha)^{2(t-j-1)}\E\norm{\zeta_j}^2\\
&+2\sum_{j=0}^{t-1}(1-\alpha)^{2t-j-1}
\underbrace{\E\inner{v_0-\nabla f(x_0)}{\zeta_j}}_{=0}\\
&+2\sum_{0\leq j<k<t}(1-\alpha)^{2t-j-k-2}
\underbrace{\E\inner{\zeta_j}{\zeta_k}}_{=0}.
\end{align*}
For the initialization, distinct samples are independent and centered, so
\begin{align*}
\E\norm{v_0-\nabla f(x_0)}^2
={}&\frac1{B^2}\sum_{k=1}^B
\E\norm{\nabla f(x_0;\xi_{0,k})-\nabla f(x_0)}^2\\
&+\frac2{B^2}\sum_{1\leq k<l\leq B}
\underbrace{\E\inner{\nabla f(x_0;\xi_{0,k})-\nabla f(x_0)}
{\nabla f(x_0;\xi_{0,l})-\nabla f(x_0)}}_{=0}\\
\leq{}&\frac{B\sigma^2}{B^2}=\frac{\sigma^2}{B}.
\end{align*}
Substituting this bound and \eqref{eq:online-noise} into the square expansion yields
\begin{align}
\E\norm{u_t-\nabla f(x_t)}^2
\leq{}&(1-\alpha)^{2t}\frac{\sigma^2}{B}\nonumber\\
&+\sum_{j=0}^{t-1}(1-\alpha)^{2(t-j-1)}
\left(2L^2\E\norm{x_{j+1}-x_j}^2+2\alpha^2\sigma^2\right).
\label{eq:virtual-pointwise}
\end{align}

We now sum this inequality over $t$. The initial error has total weight
\[
\sum_{t=0}^{T-1}(1-\alpha)^{2t}
=\frac{1-(1-\alpha)^{2T}}{\alpha(2-\alpha)}
\leq\min\left\{T,\frac1\alpha\right\}.
\]
The first bound follows because each summand is at most one. The second follows by extending the geometric sum to infinity and using $2-\alpha\geq1$.
For each subsequent innovation $j$, its total weight is at most
\[
\sum_{t=j+1}^{T-1}(1-\alpha)^{2(t-j-1)}
\leq\sum_{k=0}^{\infty}(1-\alpha)^{2k}
=\frac{1}{\alpha(2-\alpha)}\leq\frac{1}{\alpha}.
\]
Changing the order of summation uses
\[
\{(t,j):0\leq t\leq T-1,\ 0\leq j<t\}
=\{(t,j):0\leq j\leq T-2,\ j+1\leq t\leq T-1\}.
\]
Therefore,
\begin{align*}
&\sum_{t=0}^{T-1}\sum_{j=0}^{t-1}(1-\alpha)^{2(t-j-1)}
\left(2L^2\E\norm{x_{j+1}-x_j}^2+2\alpha^2\sigma^2\right)\\
&\quad=\sum_{j=0}^{T-2}
\left(2L^2\E\norm{x_{j+1}-x_j}^2+2\alpha^2\sigma^2\right)
\sum_{t=j+1}^{T-1}(1-\alpha)^{2(t-j-1)}\\
&\quad\leq\frac{2L^2}{\alpha}\sum_{j=0}^{T-2}\E\norm{x_{j+1}-x_j}^2
+2\alpha\sigma^2(T-1)\\
&\quad\leq\frac{2L^2}{\alpha}\sum_{j=0}^{T-1}\E\norm{x_{j+1}-x_j}^2
+2\alpha\sigma^2T.
\end{align*}
All sums with upper endpoint $T-2$ are empty when $T=1$. Adding the initialization term proves \eqref{eq:virtual-error}. Its contribution after division by $T$ is
\[
\frac{\sigma^2}{BT}\min\left\{T,\frac1\alpha\right\}
=\frac{\sigma^2}{B}\min\left\{1,\frac1{\alpha T}\right\}.
\]

To bound the correction itself, first expand the second moment of a stochastic gradient:
\begin{align*}
\E\norm{\nabla f(x_t;\xi_t)}^2
={}&\E\norm{\nabla f(x_t)}^2
+\E\norm{\nabla f(x_t;\xi_t)-\nabla f(x_t)}^2\\
&+2\underbrace{\E\inner{\nabla f(x_t)}
{\nabla f(x_t;\xi_t)-\nabla f(x_t)}}_{=0}\\
\leq{}&\E\norm{\nabla f(x_t)}^2+\sigma^2.
\end{align*}
The second expression in \eqref{eq:online-correction} then gives
\begin{align*}
\E\norm{\delta_t}^2
&\leq2\E\norm{\nabla f(x_{t+1};\xi_t)-\nabla f(x_t;\xi_t)}^2
+2\alpha^2\E\norm{\nabla f(x_t;\xi_t)}^2\\
&\leq2L^2\E\norm{x_{t+1}-x_t}^2
+2\alpha^2\left(\E\norm{\nabla f(x_t)}^2+\sigma^2\right).
\end{align*}
This is \eqref{eq:online-packet}.
\end{proof}

\begin{lemma}
\label{lem:online-error}
Under Assumption~\ref{ass:stochastic}, the sampling protocol of Section~\ref{sec:setting}, and the delay bound \eqref{eq:delay}, Algorithm~\ref{alg:online} with any $0<\alpha<1$ satisfies
\begin{align}
\sum_{t=0}^{T-1}\E\norm{v_t-\nabla f(x_t)}^2
\leq{}&\frac{2\sigma^2}{B}\min\left\{T,\frac1\alpha\right\}
+4L^2\left(\frac1\alpha+\tau^2\right)
\sum_{t=0}^{T-1}\E\norm{x_{t+1}-x_t}^2\nonumber\\
&+4\alpha^2\tau^2\sum_{t=0}^{T-1}\E\norm{\nabla f(x_t)}^2
+4\sigma^2(\alpha+\alpha^2\tau^2)T.
\label{eq:online-error}
\end{align}
\end{lemma}

\begin{proof}
For $j\in\mathcal M_t$, the integer arrival time satisfies
\[
t+1\leq A_j\leq j+1+\tau,
\qquad\text{hence}\qquad j\geq t-\tau.
\]
Together with $0\leq j<t$, this gives
\[
\mathcal M_t\subseteq\{\max\{0,t-\tau\},\ldots,t-1\},
\qquad |\mathcal M_t|\leq\tau.
\]
If $\tau=0$, all missing sets are empty. Otherwise, Cauchy--Schwarz and $0\leq(1-\alpha)^{t-j-1}\leq1$ give, for each realized execution,
\begin{align*}
\norm{\sum_{j\in\mathcal M_t}(1-\alpha)^{t-j-1}\delta_j}^2
&\leq|\mathcal M_t|\sum_{j\in\mathcal M_t}(1-\alpha)^{2(t-j-1)}\norm{\delta_j}^2\\
&\leq\tau\sum_{j=\max\{0,t-\tau\}}^{t-1}\norm{\delta_j}^2.
\end{align*}
The index condition in the last sum can also be written as
\[
\max\{0,t-\tau\}\leq j\leq t-1
\quad\Longleftrightarrow\quad
0\leq j\leq T-2,\quad j+1\leq t\leq\min\{T-1,j+\tau\},
\]
when $0\leq t\leq T-1$. Thus, after taking expectations,
\begin{align*}
\sum_{t=0}^{T-1}\sum_{j=\max\{0,t-\tau\}}^{t-1}\E\norm{\delta_j}^2
&=\sum_{j=0}^{T-2}\sum_{t=j+1}^{\min\{T-1,j+\tau\}}\E\norm{\delta_j}^2\\
&\leq\tau\sum_{j=0}^{T-2}\E\norm{\delta_j}^2
\leq\tau\sum_{j=0}^{T-1}\E\norm{\delta_j}^2.
\end{align*}
The inner sum contains at most $\tau$ terms. Combining these inequalities with \eqref{eq:online-packet} yields
\begin{align}
&\sum_{t=0}^{T-1}\E\norm{\sum_{j\in\mathcal M_t}(1-\alpha)^{t-j-1}\delta_j}^2\nonumber\\
&\quad\leq\tau^2\sum_{j=0}^{T-1}\E\norm{\delta_j}^2\nonumber\\
&\quad\leq2L^2\tau^2\sum_{j=0}^{T-1}\E\norm{x_{j+1}-x_j}^2
+2\alpha^2\tau^2\sum_{j=0}^{T-1}\E\norm{\nabla f(x_j)}^2
+2\alpha^2\tau^2\sigma^2T.
\label{eq:online-missing}
\end{align}
Finally, \eqref{eq:online-decomposition} and $\norm{a-b}^2\leq2\norm{a}^2+2\norm{b}^2$ imply
\begin{align*}
&\sum_{t=0}^{T-1}\E\norm{v_t-\nabla f(x_t)}^2\\
&\quad\leq2\sum_{t=0}^{T-1}\E\norm{u_t-\nabla f(x_t)}^2
+2\sum_{t=0}^{T-1}\E\norm{\sum_{j\in\mathcal M_t}(1-\alpha)^{t-j-1}\delta_j}^2\\
&\quad\leq\frac{2\sigma^2}{B}\min\left\{T,\frac1\alpha\right\}
+\left(\frac{4L^2}{\alpha}+4L^2\tau^2\right)
\sum_{t=0}^{T-1}\E\norm{x_{t+1}-x_t}^2\\
&\qquad\quad+4\alpha^2\tau^2\sum_{t=0}^{T-1}\E\norm{\nabla f(x_t)}^2
+\bigl(4\alpha\sigma^2+4\alpha^2\tau^2\sigma^2\bigr)T,
\end{align*}
which is \eqref{eq:online-error}. Arrival indicators were removed by a bound valid for every execution, before taking expectations.
\end{proof}

\subsection{Proof of Theorem~\ref{thm:online}}

By \eqref{eq:update},
\[
\sum_{t=0}^{T-1}\E\norm{x_{t+1}-x_t}^2
=\eta^2\sum_{t=0}^{T-1}\E\norm{v_t}^2.
\]
Substituting \eqref{eq:online-error} into \eqref{eq:descent-sum} therefore gives
\begin{align*}
\sum_{t=0}^{T-1}\E\norm{\nabla f(x_t)}^2
\leq{}&\frac{2\Delta_f}{\eta}
-(1-L\eta)\sum_{t=0}^{T-1}\E\norm{v_t}^2
+\frac{2\sigma^2}{B}\min\left\{T,\frac1\alpha\right\}\\
&+4L^2\eta^2\left(\frac1\alpha+\tau^2\right)
\sum_{t=0}^{T-1}\E\norm{v_t}^2\\
&+4\alpha^2\tau^2\sum_{t=0}^{T-1}\E\norm{\nabla f(x_t)}^2
+4\sigma^2(\alpha+\alpha^2\tau^2)T.
\end{align*}
Move the gradient and estimator sums on the right to the left:
\begin{align}
&(1-4\alpha^2\tau^2)\sum_{t=0}^{T-1}\E\norm{\nabla f(x_t)}^2\nonumber\\
&\quad+\left[1-L\eta-4L^2\eta^2\left(\frac1\alpha+\tau^2\right)\right]
\sum_{t=0}^{T-1}\E\norm{v_t}^2\nonumber\\
&\qquad\leq\frac{2\Delta_f}{\eta}
+\frac{2\sigma^2}{B}\min\left\{T,\frac1\alpha\right\}
+4\sigma^2(\alpha+\alpha^2\tau^2)T.
\label{eq:online-master}
\end{align}
We verify the coefficients and the remaining terms separately.

For the gradient coefficient, \eqref{eq:online-params-a} gives
\[
4\alpha^2\tau^2
=\frac{\tau^2}{64(1+\sigma/\epsilon)^2(1+\sigma/\epsilon+\tau)^2}
\leq\frac1{64}.
\]
For the estimator coefficient, the same parameters give
\begin{align*}
L\eta&=\frac1{16(1+\sigma/\epsilon+\tau)}\leq\frac1{16},\\
\frac{4L^2\eta^2}{\alpha}
&=\frac{1+\sigma/\epsilon}{4(1+\sigma/\epsilon+\tau)}\leq\frac14,\\
4L^2\eta^2\tau^2
&=\frac{\tau^2}{64(1+\sigma/\epsilon+\tau)^2}\leq\frac1{64}.
\end{align*}
It follows that
\begin{align*}
1-4\alpha^2\tau^2&\geq\frac{63}{64},\\
1-L\eta-4L^2\eta^2\left(\frac1\alpha+\tau^2\right)
&\geq1-\frac1{16}-\frac14-\frac1{64}=\frac{43}{64}>0.
\end{align*}
The estimator sum on the left of \eqref{eq:online-master} can therefore be dropped.

The choices of $B$ and $T$ imply
\begin{align*}
B&\geq\frac{16\sigma^2}{\epsilon^2}
\min\left\{1,\frac1{\alpha T}\right\},\\
\frac{2\sigma^2}{B}\min\left\{1,\frac1{\alpha T}\right\}
&\leq\frac{\epsilon^2}{8},\\
T&\geq\frac{64L\Delta_f(1+\sigma/\epsilon+\tau)}{\epsilon^2},\\
\frac{2\Delta_f}{\eta T}
&=\frac{32L\Delta_f(1+\sigma/\epsilon+\tau)}{T}
\leq\frac{\epsilon^2}{2}.
\end{align*}
When $\sigma=0$ or $\Delta_f=0$, the corresponding numerator is zero, so these bounds still hold. For the remaining terms, use
\[
\frac{\sigma}{1+\sigma/\epsilon}\leq\epsilon,
\qquad
1+\frac\sigma\epsilon\leq1+\frac\sigma\epsilon+\tau,
\qquad
\frac{\tau}{1+\sigma/\epsilon+\tau}\leq1.
\]
Then
\begin{align*}
4\alpha\sigma^2
&=\frac{\sigma^2}{4(1+\sigma/\epsilon)(1+\sigma/\epsilon+\tau)}
\leq\frac{\sigma^2}{4(1+\sigma/\epsilon)^2}
\leq\frac{\epsilon^2}{4},\\
4\alpha^2\tau^2\sigma^2
&=\frac{\sigma^2\tau^2}{64(1+\sigma/\epsilon)^2(1+\sigma/\epsilon+\tau)^2}
\leq\frac{\epsilon^2}{64},
\end{align*}
Dividing \eqref{eq:online-master} by $T(1-4\alpha^2\tau^2)$ gives
\[
\frac1T\sum_{t=0}^{T-1}\E\norm{\nabla f(x_t)}^2
\leq
\frac{\epsilon^2/2+\epsilon^2/8+\epsilon^2/4+\epsilon^2/64}{63/64}
=\frac{19}{21}\epsilon^2.
\]
Since the output index is uniform and independent,
\[
\E\norm{\nabla f(x_R)}^2
=\sum_{t=0}^{T-1}\frac1T\E\norm{\nabla f(x_t)}^2
\leq\frac{19}{21}\epsilon^2\leq\epsilon^2.
\]
For the cost, initialization issues $B$ single-point evaluations. Each of the $T$ updates issues one two-point task, including the last update. The ceiling bounds give
\begin{align*}
B&\leq1+\frac{16\sigma^2}{\epsilon^2},\\
T&\leq1+\frac{64L\Delta_f(1+\sigma/\epsilon+\tau)}{\epsilon^2},\\
B+2T&\leq3+\frac{16\sigma^2}{\epsilon^2}
+\frac{128L\Delta_f(1+\sigma/\epsilon+\tau)}{\epsilon^2}.
\end{align*}
This count includes every issued evaluation, even if its response arrives after the algorithm returns, and proves \eqref{eq:online-complexity}. \qed

\section{Variance Reduction for Finite-Sum Optimization}
\label{sec:finite}

We next consider the finite-sum objective
\begin{equation}
f(x)=\frac{1}{n}\sum_{i=1}^n f_i(x).
\label{eq:finite-objective}
\end{equation}
Here an exact full gradient can be computed with $n$ component-gradient evaluations. We use periodic exact refreshes to control the accumulated estimation error, while continuing to incorporate delayed corrections between refreshes.

\begin{assumption}[Finite-sum average smoothness]
\label{ass:finite-smooth}
For all $x,y\in\R^d$, and $L>0$
\[
\frac{1}{n}\sum_{i=1}^n\norm{\nabla f_i(x)-\nabla f_i(y)}^2
\leq L^2\norm{x-y}^2,
\]
\end{assumption}

By Jensen's inequality, Assumption~\ref{ass:finite-smooth} implies that $f$ is $L$-smooth. Individual $L$-smoothness of the components is sufficient, but is not required. We do not assume bounded component gradients or a uniform bound on their variance.

\subsection{The proposed method}

Let $q$ be the refresh period. At rounds $0,q,2q,\ldots$, we fix the current parameter, compute its exact full gradient, and discard corrections from earlier epochs. Between refreshes, task $t$ draws $b$ independent uniform indices $i_{t,1},\ldots,i_{t,b}$ and returns
\begin{equation}
\delta_t=\frac{1}{b}\sum_{k=1}^b
\bigl(\nabla f_{i_{t,k}}(x_{t+1})-\nabla f_{i_{t,k}}(x_t)\bigr).
\label{eq:finite-correction}
\end{equation}
For a given round $t$, let $r_t=q\lfloor t/q\rfloor$ be the first round of its epoch. The implemented estimator is
\begin{equation}
v_t=\nabla f(x_{r_t})+
\sum_{j=r_t}^{t-1}\ind_{\{A_j\leq t\}}\delta_j.
\label{eq:finite-estimator}
\end{equation}
If all corrections have returned, their population counterparts telescope from $x_{r_t}$ to $x_t$. With delayed feedback, we retain the available corrections and continue updating. Arrival order does not affect the sum.

\begin{algorithm}[H]
\caption{Delayed recursive variance reduction for finite sums}
\label{alg:finite}
\begin{algorithmic}[1]
\Require Initial point $x_0$, update number $T$, period $q$, batch size $b$, learning rate $\eta$
\For{$t=0,\ldots,T-1$}
  \If{$t$ is a multiple of $q$}
    \State Fix $x_t$, compute $\nabla f(x_t)$, and set $v_t=\nabla f(x_t)$.
    \State Discard all corrections from preceding epochs.
  \Else
    \State $v_t=v_{t-1}+\displaystyle\sum_{j\geq r_t:A_j=t}\delta_j$.
  \EndIf
  \State Update $x_{t+1}=x_t-\eta v_t$.
  \State Draw $b$ fresh uniform indices and issue task \eqref{eq:finite-correction}.
\EndFor
\State Choose $R$ uniformly from $\{0,\ldots,T-1\}$ and return $x_R$.
\end{algorithmic}
\end{algorithm}

\subsection{Convergence guarantee}

\begin{theorem}
\label{thm:finite}
Under Assumptions~\ref{ass:lower} and~\ref{ass:finite-smooth}, suppose that \eqref{eq:delay} holds. For any positive integers $q,b,T$, Algorithm~\ref{alg:finite} with
\begin{equation}
\eta=\frac{1}{2L(1+\sqrt{q/b}+\min\{\tau,q-1\})}
\label{eq:finite-eta}
\end{equation}
guarantees
\begin{equation}
\E\norm{\nabla f(x_R)}^2
=\frac{1}{T}\sum_{t=0}^{T-1}\E\norm{\nabla f(x_t)}^2
\leq\frac{2\Delta_f}{\eta T}.
\label{eq:finite-rate}
\end{equation}
The number of issued component-gradient evaluations is at most
\begin{equation}
n\left\lceil\frac{T}{q}\right\rceil+2bT.
\label{eq:finite-cost}
\end{equation}
\end{theorem}

\begin{corollary}
\label{cor:finite}
Set $q=n$, $b=1$, choose $\eta$ as in \eqref{eq:finite-eta}, and let
\[
T=\max\left\{1,\left\lceil\frac{2\Delta_f}{\eta\epsilon^2}\right\rceil\right\}.
\]
Then Algorithm~\ref{alg:finite} satisfies \eqref{eq:target} with at most $n+3T$ evaluations. Its oracle complexity is
\begin{equation}
O\!\left(n+\frac{L\Delta_f(\sqrt n+\min\{\tau,n-1\})}{\epsilon^2}\right).
\label{eq:finite-complexity}
\end{equation}
\end{corollary}

\begin{remark}
For $\tau=O(\sqrt n)$, the order in \eqref{eq:finite-complexity} agrees with the usual finite-sum variance-reduction complexity. The delay contribution is at most $n-1$, because a correction ceases to affect the estimator at the next exact refresh.
\end{remark}

\begin{corollary}
\label{cor:finite-batch}
For an integer $1\leq b\leq\sqrt n$, set $q=\lceil n/b\rceil$, choose $\eta$ as in \eqref{eq:finite-eta}, and let
$T=\max\{1,\lceil2\Delta_f/(\eta\epsilon^2)\rceil\}$.
The number of updates and the number $Q$ of issued evaluations satisfy
\begin{align}
T&=O\!\left(1+\frac{L\Delta_f}{\epsilon^2}
\left(1+\frac{\sqrt n}{b}+\tau\right)\right),
\label{eq:finite-batch-updates}\\
Q&=O\!\left(n+\frac{L\Delta_f}{\epsilon^2}
\left(\sqrt n+b(1+\tau)\right)\right).
\label{eq:finite-batch-work}
\end{align}
\end{corollary}

Thus any choice $b\leq\sqrt n/(1+\tau)$ preserves the order
$O(n+L\Delta_f\sqrt n/\epsilon^2)$ of the oracle work while reducing the bound on the number of updates. Turning this tradeoff into an elapsed-time guarantee requires a service model for the batches and refreshes.

\subsection{Proof of Theorem~\ref{thm:finite}}

We first control the estimation error. For round $t$, denote the unfinished corrections from its current epoch by
\begin{equation}
\mathcal M_t=\{j:r_t\leq j<t,\ A_j>t\}.
\label{eq:finite-pending}
\end{equation}
Under \eqref{eq:delay}, every index in $\mathcal M_t$ lies between $\max\{r_t,t-\tau\}$ and $t-1$. There are at most $q-1$ preceding rounds in the current epoch, so $|\mathcal M_t|\leq\min\{\tau,q-1\}$.

\begin{lemma}
\label{lem:finite-error}
Under the conditions of Theorem~\ref{thm:finite},
\begin{align}
\sum_{t=0}^{T-1}\E\norm{v_t-\nabla f(x_t)}^2
\leq 2L^2\left(\frac{q-1}{b}+\min\{\tau,q-1\}^2\right)
\sum_{t=0}^{T-1}\E\norm{x_{t+1}-x_t}^2.
\label{eq:finite-error}
\end{align}
\end{lemma}

\begin{proof}
For task $j$, define the centered batch error by
\[
\zeta_j=\delta_j-
\bigl(\nabla f(x_{j+1})-\nabla f(x_j)\bigr).
\]
Both query points are fixed before the indices of this task are drawn. For any fixed $x,y$, uniform sampling gives
\[
\frac1n\sum_{i=1}^n\bigl(\nabla f_i(y)-\nabla f_i(x)\bigr)
=\nabla f(y)-\nabla f(x).
\]
Thus every centered sampled difference has zero mean when averaging its new index. Cross terms between distinct sampled differences also have zero mean, because the indices in the batch are independent. Expanding the square of the batch average therefore gives
\begin{align*}
\E\norm{\zeta_j}^2
&=\frac1{b^2}\sum_{k=1}^b\E
\left\|\begin{aligned}
&\nabla f_{i_{j,k}}(x_{j+1})-\nabla f_{i_{j,k}}(x_j)\\[-2pt]
&\quad-\bigl(\nabla f(x_{j+1})-\nabla f(x_j)\bigr)
\end{aligned}\right\|^2\\
&=\frac1{bn}\sum_{i=1}^n
\E\norm{\nabla f_i(x_{j+1})-\nabla f_i(x_j)}^2\\
&\quad-\frac1b\E\norm{\nabla f(x_{j+1})-\nabla f(x_j)}^2\\
&\leq\frac1{bn}\sum_{i=1}^n
\E\norm{\nabla f_i(x_{j+1})-\nabla f_i(x_j)}^2.
\end{align*}
The second equality uses uniform sampling and the identity
$\E\norm{Z-\E Z}^2=\E\norm Z^2-\norm{\E Z}^2$ with the query points held fixed. Assumption~\ref{ass:finite-smooth} now gives
\begin{equation}
\E\norm{\zeta_j}^2
\leq\frac{L^2}{b}\E\norm{x_{j+1}-x_j}^2.
\label{eq:finite-noise}
\end{equation}
We also need a bound on $\delta_j$, which includes its mean. Convexity of the squared norm and uniform sampling give
\begin{align}
\E\norm{\delta_j}^2
&\leq\frac1b\sum_{k=1}^b
\E\norm{\nabla f_{i_{j,k}}(x_{j+1})-\nabla f_{i_{j,k}}(x_j)}^2\notag\\
&=\frac1n\sum_{i=1}^n
\E\norm{\nabla f_i(x_{j+1})-\nabla f_i(x_j)}^2\notag\\
&\leq L^2\E\norm{x_{j+1}-x_j}^2.
\label{eq:finite-packet}
\end{align}

The full sequence of population differences telescopes:
\[
\nabla f(x_t)=\nabla f(x_{r_t})+
\sum_{j=r_t}^{t-1}\bigl(\nabla f(x_{j+1})-\nabla f(x_j)\bigr).
\]
To subtract this identity from \eqref{eq:finite-estimator}, first write the sum of returned corrections as
\[
\sum_{\substack{r_t\leq j<t\\A_j\leq t}}\delta_j
=\sum_{j=r_t}^{t-1}\delta_j-\sum_{j\in\mathcal M_t}\delta_j.
\]
It follows that
\begin{equation}
v_t-\nabla f(x_t)
=\sum_{j=r_t}^{t-1}\zeta_j-
\sum_{j\in\mathcal M_t}\delta_j.
\label{eq:finite-decomposition}
\end{equation}
The first sum contains all generated sampling errors in the epoch, rather than just the returned errors. This allows us to control it using fresh sampling even when completion times depend on the samples.

For $i<j$, the vector $\zeta_i$ and the point pair of task $j$ have already been determined when the new indices of task $j$ are drawn. Averaging these new indices first gives
$\E\inner{\zeta_i}{\zeta_j}=0$. This argument does not require the two point pairs to be independent. Since $r_t$ is a deterministic epoch boundary, expansion of the squared norm gives
\begin{align}
\E\norm{\sum_{j=r_t}^{t-1}\zeta_j}^2
&=\sum_{j=r_t}^{t-1}\E\norm{\zeta_j}^2
+2\sum_{r_t\leq i<j<t}\E\inner{\zeta_i}{\zeta_j}\notag\\
&=\sum_{j=r_t}^{t-1}\E\norm{\zeta_j}^2\notag\\
&\leq\frac{L^2}{b}\sum_{j=r_t}^{t-1}\E\norm{x_{j+1}-x_j}^2.
\label{eq:finite-complete-error}
\end{align}
To sum this inequality over $t$, fix an index $j$. It appears only after it is generated and before the next refresh, namely,
\[
j+1\leq t\leq q(\lfloor j/q\rfloor+1)-1.
\]
There are $q(\lfloor j/q\rfloor+1)-j-1\leq q-1$ such rounds. In particular, the correction issued in the last round of an epoch appears zero times. Reversing the order of summation makes this count explicit:
\begin{align*}
&\sum_{t=0}^{T-1}\sum_{j=r_t}^{t-1}
\E\norm{x_{j+1}-x_j}^2\\
&\quad=\sum_{j=0}^{T-1}
\sum_{t=j+1}^{\min\{T-1,q(\lfloor j/q\rfloor+1)-1\}}
\E\norm{x_{j+1}-x_j}^2\\
&\quad\leq(q-1)\sum_{j=0}^{T-1}
\E\norm{x_{j+1}-x_j}^2.
\end{align*}
Here and below, a sum is zero if its upper endpoint is smaller than its lower endpoint. Applying this bound to \eqref{eq:finite-complete-error} yields
\begin{equation}
\sum_{t=0}^{T-1}\E\norm{\sum_{j=r_t}^{t-1}\zeta_j}^2
\leq\frac{L^2(q-1)}{b}\sum_{j=0}^{T-1}\E\norm{x_{j+1}-x_j}^2.
\label{eq:finite-complete-sum}
\end{equation}

To bound the second sum in \eqref{eq:finite-decomposition}, we work with each realized execution before taking expectations. Cauchy--Schwarz gives
\begin{align}
\norm{\sum_{j\in\mathcal M_t}\delta_j}^2
&\leq |\mathcal M_t|\sum_{j\in\mathcal M_t}\norm{\delta_j}^2\notag\\
&\leq\min\{\tau,q-1\}\sum_{j=\max\{r_t,t-\tau\}}^{t-1}\norm{\delta_j}^2.
\label{eq:finite-missing-point}
\end{align}
The last sum has deterministic endpoints, so taking expectations does not retain a sample-dependent arrival indicator. For a fixed $j$, membership in this interval requires both $j<t\leq j+\tau$ and that $t$ remain in the epoch of $j$. Consequently,
\begin{align*}
&\sum_{t=0}^{T-1}\sum_{j=\max\{r_t,t-\tau\}}^{t-1}
\E\norm{\delta_j}^2\\
&\quad=\sum_{j=0}^{T-1}
\sum_{t=j+1}^{\min\{T-1,j+\tau,q(\lfloor j/q\rfloor+1)-1\}}
\E\norm{\delta_j}^2\\
&\quad\leq\min\{\tau,q-1\}\sum_{j=0}^{T-1}\E\norm{\delta_j}^2.
\end{align*}
The inner sum contains at most $\tau$ rounds because $t\leq j+\tau$, and at most $q-1$ rounds because of the next refresh. Summing \eqref{eq:finite-missing-point} and applying this count gives
\begin{align}
\sum_{t=0}^{T-1}\E\norm{\sum_{j\in\mathcal M_t}\delta_j}^2
&\leq\min\{\tau,q-1\}^2\sum_{j=0}^{T-1}\E\norm{\delta_j}^2
\notag\\
&\leq L^2\min\{\tau,q-1\}^2\sum_{j=0}^{T-1}\E\norm{x_{j+1}-x_j}^2.
\label{eq:finite-missing-sum}
\end{align}
Empty sums are zero, so these inequalities also cover $\tau=0$ and $q=1$.
Finally, apply $\norm{a-b}^2\leq2\norm{a}^2+2\norm{b}^2$ to \eqref{eq:finite-decomposition} and sum:
\begin{align*}
\sum_{t=0}^{T-1}\E\norm{v_t-\nabla f(x_t)}^2
&\leq2\sum_{t=0}^{T-1}\E\norm{\sum_{j=r_t}^{t-1}\zeta_j}^2
+2\sum_{t=0}^{T-1}\E\norm{\sum_{j\in\mathcal M_t}\delta_j}^2\\
&\leq2L^2\left(\frac{q-1}{b}+\min\{\tau,q-1\}^2\right)
\sum_{j=0}^{T-1}\E\norm{x_{j+1}-x_j}^2.
\end{align*}
This is \eqref{eq:finite-error}.
\end{proof}

\begin{proof}[Proof of Theorem~\ref{thm:finite}]
The update gives $\norm{x_{t+1}-x_t}^2=\eta^2\norm{v_t}^2$. We first use \eqref{eq:descent-sum} to separate the estimation error from the descent term:
\begin{align*}
\sum_{t=0}^{T-1}\E\norm{\nabla f(x_t)}^2
&\leq\frac{2\Delta_f}{\eta}
-(1-L\eta)\sum_{t=0}^{T-1}\E\norm{v_t}^2
+\sum_{t=0}^{T-1}\E\norm{v_t-\nabla f(x_t)}^2\\
&\leq\frac{2\Delta_f}{\eta}
-(1-L\eta)\sum_{t=0}^{T-1}\E\norm{v_t}^2\\
&\quad+2L^2\eta^2\left(\frac{q-1}{b}+\min\{\tau,q-1\}^2\right)
\sum_{t=0}^{T-1}\E\norm{v_t}^2,
\end{align*}
where the second inequality is Lemma~\ref{lem:finite-error}. Collecting the two coefficients of $\sum_t\E\norm{v_t}^2$ yields
\begin{align}
\sum_{t=0}^{T-1}\E\norm{\nabla f(x_t)}^2
\leq{}&\frac{2\Delta_f}{\eta}
-\left[1-L\eta-2L^2\eta^2
\left(\frac{q-1}{b}+\min\{\tau,q-1\}^2\right)\right]
\sum_{t=0}^{T-1}\E\norm{v_t}^2.
\label{eq:finite-before-eta}
\end{align}
For the learning rate in \eqref{eq:finite-eta},
\[
L\eta=\frac{1}{2(1+\sqrt{q/b}+\min\{\tau,q-1\})}\leq\frac12.
\]
Moreover,
\begin{align*}
2L^2\eta^2\left(\frac{q-1}{b}+\min\{\tau,q-1\}^2\right)
&\leq\frac{q/b+\min\{\tau,q-1\}^2}
{2(1+\sqrt{q/b}+\min\{\tau,q-1\})^2}\\
&\leq\frac{(\sqrt{q/b}+\min\{\tau,q-1\})^2}
{2(1+\sqrt{q/b}+\min\{\tau,q-1\})^2}\\
&\leq\frac12.
\end{align*}
Hence the bracket in \eqref{eq:finite-before-eta} is at least $1-1/2-1/2=0$. Dropping its nonpositive contribution and dividing by $T$ gives
\[
\frac1T\sum_{t=0}^{T-1}\E\norm{\nabla f(x_t)}^2
\leq\frac{2\Delta_f}{\eta T}.
\]
Since $R$ is sampled uniformly after the run, its output satisfies
\[
\E\norm{\nabla f(x_R)}^2
=\sum_{t=0}^{T-1}\frac1T\E\norm{\nabla f(x_t)}^2,
\]
which proves \eqref{eq:finite-rate}.

The full gradients are requested at $0,q,\ldots,q\lfloor(T-1)/q\rfloor$. Their number is $1+\lfloor(T-1)/q\rfloor=\lceil T/q\rceil$, and each uses $n$ evaluations. Each of the $T$ rounds also issues one task containing $b$ gradient pairs. Counting both gradients in every pair, including pairs discarded after a refresh or still unfinished at termination, gives $2bT$ evaluations. Adding the two costs proves \eqref{eq:finite-cost}.
\end{proof}

\begin{proof}[Proof of Corollary~\ref{cor:finite}]
The definition of $T$ ensures $T\geq2\Delta_f/(\eta\epsilon^2)$, so Theorem~\ref{thm:finite} gives the required accuracy. It also ensures
\[
T\leq1+\frac{2\Delta_f}{\eta\epsilon^2}
=1+\frac{4L\Delta_f}{\epsilon^2}
\bigl(1+\sqrt n+\min\{\tau,n-1\}\bigr).
\]
This bound includes $\Delta_f=0$, in which case $T=1$. For the oracle count,
\begin{align*}
n\left\lceil\frac Tn\right\rceil+2T
&\leq n\left(\frac Tn+1\right)+2T=n+3T\\
&\leq n+3+\frac{12L\Delta_f}{\epsilon^2}
\bigl(1+\sqrt n+\min\{\tau,n-1\}\bigr).
\end{align*}
Using $n+3\leq4n$ and $1+\sqrt n\leq2\sqrt n$ gives \eqref{eq:finite-complexity}.
\end{proof}

\begin{proof}[Proof of Corollary~\ref{cor:finite-batch}]
Since $q=\lceil n/b\rceil$, we have
\[
\frac nb\leq q\leq\frac nb+1,
\qquad
\sqrt{\frac qb}
\leq\sqrt{\frac n{b^2}+\frac1b}
\leq\frac{\sqrt n}{b}+\frac1{\sqrt b}
\leq\frac{\sqrt n}{b}+1.
\]
Also, $\min\{\tau,q-1\}\leq\tau$. Substitution into \eqref{eq:finite-eta} therefore gives
\begin{align*}
T&\leq1+\frac{2\Delta_f}{\eta\epsilon^2}\\
&=1+\frac{4L\Delta_f}{\epsilon^2}
\bigl(1+\sqrt{q/b}+\min\{\tau,q-1\}\bigr)\\
&\leq1+\frac{4L\Delta_f}{\epsilon^2}
\left(2+\frac{\sqrt n}{b}+\tau\right),
\end{align*}
which proves \eqref{eq:finite-batch-updates}. For the evaluation count, $q\geq n/b$ implies $n/q\leq b$, and hence
\begin{align*}
Q&\leq n\left\lceil\frac Tq\right\rceil+2bT\\
&\leq n+\left(\frac nq+2b\right)T\leq n+3bT\\
&\leq n+3b+\frac{12L\Delta_f}{\epsilon^2}
\bigl(2b+\sqrt n+b\tau\bigr).
\end{align*}
Finally, $b\leq\sqrt n\leq n$ gives $n+3b\leq4n$, proving \eqref{eq:finite-batch-work}. The target accuracy follows from the same choice $T\geq2\Delta_f/(\eta\epsilon^2)$ and Theorem~\ref{thm:finite}.
\end{proof}

\section{Finite-Sum Updates without Refresh Waits}
\label{sec:snapshots}

Algorithm~\ref{alg:finite} waits for an exact gradient at every refresh. We now allow these gradients to be computed while the server continues updating. The correction tasks remain unchanged. The additional requirement concerns the completion time of an exact full gradient.

\begin{assumption}
\label{ass:snapshot-service}
Fix an integer $D\geq1$. At every round $r=q,2q,\ldots$, the server can request the exact gradient $\nabla f(x_r)$ at the fixed point $x_r$. All $n$ component evaluations in this request are completed and their average is available before the update at round $r+D$. The server can continue to issue correction tasks and update the parameter while the request is being processed.
\end{assumption}

Assumption~\ref{ass:snapshot-service} is separate from \eqref{eq:delay}. The latter bounds the delay of a two-point correction and it does not bound the time needed to complete $n$ component evaluations.

\subsection{The proposed method}

We call an exact gradient at a fixed earlier point a snapshot. The snapshot requested at round $r$ is used starting at the prescribed round $r+D$, even if it returns sooner. Consequently, the snapshot used at round $t$ has index
\begin{equation}
r_t=\max\left\{0,\ q\left\lfloor\frac{t-D}{q}\right\rfloor\right\}.
\label{eq:snapshot-anchor}
\end{equation}
The value $r_t=0$ refers to the completed initialization gradient. Define
\begin{equation}
v_t=\nabla f(x_{r_t})+
\sum_{\substack{r_t\leq j<t\\A_j\leq t}}\delta_j,
\label{eq:snapshot-estimator}
\end{equation}
where $\delta_j$ is the same batch of component-gradient differences as in \eqref{eq:finite-correction}. The point continues to move while a snapshot is being computed. When that snapshot is used, the returned corrections generated since its point was fixed are retained.

\begin{algorithm}[tb]
\caption{Delayed recursive variance reduction with overlapping snapshots}
\label{alg:snapshot}
\begin{algorithmic}[1]
\Require Initial point $x_0$, update number $T$, period $q$, batch size $b$, full-gradient deadline $D\geq1$, learning rate $\eta$
\State Compute $\nabla f(x_0)$ and wait for all $n$ component evaluations.
\For{$t=0,\ldots,T-1$}
  \State Receive all feedback available at the beginning of round $t$.
  \If{$t\geq q$ and $t$ is a multiple of $q$}
    \State Request $\nabla f(x_t)$ at the fixed point $x_t$, for use at round $t+D$.
  \EndIf
  \State Set $r_t$ by \eqref{eq:snapshot-anchor} and form $v_t$ by \eqref{eq:snapshot-estimator}.
  \State Update $x_{t+1}=x_t-\eta v_t$.
  \State Draw $b$ fresh uniform indices and issue task \eqref{eq:finite-correction}.
\EndFor
\State Choose $R$ uniformly from $\{0,\ldots,T-1\}$ and return $x_R$.
\end{algorithmic}
\end{algorithm}

\subsection{Convergence guarantee}

\begin{theorem}
\label{thm:snapshot}
Under Assumptions~\ref{ass:lower}, \ref{ass:finite-smooth}, and \ref{ass:snapshot-service}, suppose that correction tasks satisfy \eqref{eq:delay}. For positive integers $q,b,T$ and $D\geq1$, Algorithm~\ref{alg:snapshot} with
\begin{equation}
\eta=\frac{1}{2L\left(1+\sqrt{(q+D)/b}+\tau\right)}
\label{eq:snapshot-eta}
\end{equation}
satisfies
\begin{equation}
\E\norm{\nabla f(x_R)}^2
\leq\frac{2\Delta_f}{\eta T}.
\label{eq:snapshot-rate}
\end{equation}
The number of issued component-gradient evaluations is at most
\begin{equation}
n\left\lceil\frac{T}{q}\right\rceil+2bT.
\label{eq:snapshot-cost}
\end{equation}
\end{theorem}

\begin{corollary}
\label{cor:snapshot}
Set $q=n$, $b=1$, choose $\eta$ as in \eqref{eq:snapshot-eta}, and let
$T=\max\{1,\lceil2\Delta_f/(\eta\epsilon^2)\rceil\}$.
Then Algorithm~\ref{alg:snapshot} satisfies \eqref{eq:target} with oracle complexity
\begin{equation}
O\!\left(n+\frac{L\Delta_f}{\epsilon^2}
\left(\sqrt n+\sqrt D+\tau\right)\right).
\label{eq:snapshot-complexity}
\end{equation}
\end{corollary}

The server does not wait for a refresh after initialization. If the service model gives $D=O(n+\tau^2)$, the order in \eqref{eq:snapshot-complexity} is $O(n+L\Delta_f(\sqrt n+\tau)/\epsilon^2)$.

\subsection{Proof of Theorem~\ref{thm:snapshot}}

The new ingredient is the time for which a sampling error remains in the estimator. The exact snapshot cancels all corrections preceding its anchor, even though that anchor is older than the current parameter.

\begin{lemma}
\label{lem:snapshot-error}
Under the conditions of Theorem~\ref{thm:snapshot},
\begin{equation}
\sum_{t=0}^{T-1}\E\norm{v_t-\nabla f(x_t)}^2
\leq 2L^2\left(\frac{q+D-1}{b}+\tau^2\right)
\sum_{t=0}^{T-1}\E\norm{x_{t+1}-x_t}^2.
\label{eq:snapshot-error}
\end{equation}
\end{lemma}

\begin{proof}
For each task $j$, define
\[
\zeta_j=\delta_j-
\bigl(\nabla f(x_{j+1})-\nabla f(x_j)\bigr),
\qquad
\mathcal M_t=\{j:r_t\leq j<t,\ A_j>t\}.
\]
The samples for task $j$ are drawn after its point pair has been fixed. Therefore, the proof of \eqref{eq:finite-noise} gives
\begin{equation}
\E\norm{\zeta_j}^2
\leq\frac{L^2}{b}\E\norm{x_{j+1}-x_j}^2,
\qquad
\E\norm{\delta_j}^2
\leq L^2\E\norm{x_{j+1}-x_j}^2.
\label{eq:snapshot-packets}
\end{equation}
The first inequality follows by expanding the squared norm of the average of $b$ centered differences; its $b(b-1)$ cross terms have zero mean. The second is \eqref{eq:finite-packet}, which bounds the uncentered correction. These bounds use only fresh sampling after the point pair is chosen, so they hold for the actual iterates of Algorithm~\ref{alg:snapshot} as well.

The population differences telescope from the snapshot point to the current point:
\[
\nabla f(x_t)=\nabla f(x_{r_t})+
\sum_{j=r_t}^{t-1}
\bigl(\nabla f(x_{j+1})-\nabla f(x_j)\bigr).
\]
The returned corrections and missing corrections partition all indices from $r_t$ to $t-1$. Therefore,
\begin{align*}
v_t-\nabla f(x_t)
&=\sum_{j=r_t}^{t-1}\delta_j-\sum_{j\in\mathcal M_t}\delta_j\\
&\quad-\sum_{j=r_t}^{t-1}
\bigl(\nabla f(x_{j+1})-\nabla f(x_j)\bigr)\\
&=\sum_{j=r_t}^{t-1}
\left[\delta_j-\bigl(\nabla f(x_{j+1})-\nabla f(x_j)\bigr)\right]
-\sum_{j\in\mathcal M_t}\delta_j.
\end{align*}
By the definition of $\zeta_j$, this is the decomposition
\begin{equation}
v_t-\nabla f(x_t)
=\sum_{j=r_t}^{t-1}\zeta_j-
\sum_{j\in\mathcal M_t}\delta_j.
\label{eq:snapshot-decomposition}
\end{equation}

For the first sum, the endpoints $r_t$ and $t-1$ are deterministic. If $i<j$, then $\zeta_i$ and the point pair of task $j$ are determined before the indices of task $j$ are drawn. Averaging those new indices first gives $\E\inner{\zeta_i}{\zeta_j}=0$. Expanding the square and using \eqref{eq:snapshot-packets}, we obtain
\begin{align}
\E\norm{\sum_{j=r_t}^{t-1}\zeta_j}^2
&=\sum_{j=r_t}^{t-1}\E\norm{\zeta_j}^2
+2\sum_{r_t\leq i<j<t}\E\inner{\zeta_i}{\zeta_j}\notag\\
&=\sum_{j=r_t}^{t-1}\E\norm{\zeta_j}^2\notag\\
&\leq\frac{L^2}{b}\sum_{j=r_t}^{t-1}
\E\norm{x_{j+1}-x_j}^2.
\label{eq:snapshot-complete-point}
\end{align}

We next count the rounds containing a fixed index $j$. It enters the sum at round $j+1$. The first snapshot with anchor strictly larger than $j$ is requested at
$q(\lfloor j/q\rfloor+1)$ and activated $D$ rounds later. The anchor never decreases, so
\[
r_t\leq j<t
\quad\Longleftrightarrow\quad
j+1\leq t\leq q(\lfloor j/q\rfloor+1)+D-1.
\]
This also holds for the initial anchor $0$: the first replacement is requested at $q$ and activated at $q+D$. Thus reversal of the two sums gives
\begin{align*}
&\sum_{t=0}^{T-1}\sum_{j=r_t}^{t-1}
\E\norm{x_{j+1}-x_j}^2\\
&\quad=\sum_{j=0}^{T-1}
\sum_{t=j+1}^{\min\{T-1,q(\lfloor j/q\rfloor+1)+D-1\}}
\E\norm{x_{j+1}-x_j}^2\\
&\quad\leq(q+D-1)\sum_{j=0}^{T-1}
\E\norm{x_{j+1}-x_j}^2.
\end{align*}
Indeed, without the finite-horizon cutoff, the inner sum has
\[
q(\lfloor j/q\rfloor+1)+D-j-1
=q+D-1-(j-q\lfloor j/q\rfloor)
\leq q+D-1
\]
terms. The cutoff at $T-1$ can only reduce this count. Summing \eqref{eq:snapshot-complete-point} now yields
\begin{equation}
\sum_{t=0}^{T-1}\E\norm{\sum_{j=r_t}^{t-1}\zeta_j}^2
\leq\frac{L^2(q+D-1)}{b}
\sum_{j=0}^{T-1}\E\norm{x_{j+1}-x_j}^2.
\label{eq:snapshot-complete-sum}
\end{equation}

For the missing corrections, \eqref{eq:delay} implies
$\mathcal M_t\subseteq\{\max\{r_t,t-\tau\},\ldots,t-1\}$ and
$|\mathcal M_t|\leq\tau$. On each realized execution, Cauchy--Schwarz gives
\begin{align*}
\norm{\sum_{j\in\mathcal M_t}\delta_j}^2
&\leq|\mathcal M_t|\sum_{j\in\mathcal M_t}\norm{\delta_j}^2\\
&\leq\tau\sum_{j=\max\{r_t,t-\tau\}}^{t-1}\norm{\delta_j}^2\\
&\leq\tau\sum_{j=\max\{0,t-\tau\}}^{t-1}\norm{\delta_j}^2.
\end{align*}
The last step only adds nonnegative terms and removes the anchor from the bound. For a fixed $j$, the remaining interval condition is $j+1\leq t\leq j+\tau$. Therefore,
\begin{align*}
\sum_{t=0}^{T-1}\sum_{j=\max\{0,t-\tau\}}^{t-1}\E\norm{\delta_j}^2
&=\sum_{j=0}^{T-1}\sum_{t=j+1}^{\min\{T-1,j+\tau\}}
\E\norm{\delta_j}^2\\
&\leq\tau\sum_{j=0}^{T-1}\E\norm{\delta_j}^2.
\end{align*}
Taking expectations in the preceding pointwise bound and using this count gives
\begin{align}
\sum_{t=0}^{T-1}\E\norm{\sum_{j\in\mathcal M_t}\delta_j}^2
&\leq\tau^2\sum_{j=0}^{T-1}\E\norm{\delta_j}^2\notag\\
&\leq L^2\tau^2\sum_{j=0}^{T-1}
\E\norm{x_{j+1}-x_j}^2.
\label{eq:snapshot-missing-sum}
\end{align}
For $\tau=0$, the missing sets are empty and the same inequality holds. Finally, apply $\norm{a-b}^2\leq2\norm{a}^2+2\norm{b}^2$ to \eqref{eq:snapshot-decomposition}:
\begin{align*}
\sum_{t=0}^{T-1}\E\norm{v_t-\nabla f(x_t)}^2
&\leq2\sum_{t=0}^{T-1}\E\norm{\sum_{j=r_t}^{t-1}\zeta_j}^2
+2\sum_{t=0}^{T-1}\E\norm{\sum_{j\in\mathcal M_t}\delta_j}^2\\
&\leq2L^2\left(\frac{q+D-1}{b}+\tau^2\right)
\sum_{j=0}^{T-1}\E\norm{x_{j+1}-x_j}^2,
\end{align*}
where the last line uses \eqref{eq:snapshot-complete-sum} and \eqref{eq:snapshot-missing-sum}. This proves the lemma.
\end{proof}

\begin{proof}[Proof of Theorem~\ref{thm:snapshot}]
The update identity $\norm{x_{t+1}-x_t}^2=\eta^2\norm{v_t}^2$ and Lemma~\ref{lem:snapshot-error} give
\[
\sum_{t=0}^{T-1}\E\norm{v_t-\nabla f(x_t)}^2
\leq2L^2\eta^2\left(\frac{q+D-1}{b}+\tau^2\right)
\sum_{t=0}^{T-1}\E\norm{v_t}^2.
\]
Substituting this bound into \eqref{eq:descent-sum}, we obtain
\begin{align*}
\sum_{t=0}^{T-1}\E\norm{\nabla f(x_t)}^2
&\leq\frac{2\Delta_f}{\eta}
-(1-L\eta)\sum_{t=0}^{T-1}\E\norm{v_t}^2\\
&\quad+2L^2\eta^2\left(\frac{q+D-1}{b}+\tau^2\right)
\sum_{t=0}^{T-1}\E\norm{v_t}^2.
\end{align*}
Combining the two coefficients of $\sum_t\E\norm{v_t}^2$ gives
\begin{align}
\sum_{t=0}^{T-1}\E\norm{\nabla f(x_t)}^2
\leq{}&\frac{2\Delta_f}{\eta}
-\left[1-L\eta-2L^2\eta^2
\left(\frac{q+D-1}{b}+\tau^2\right)\right]
\sum_{t=0}^{T-1}\E\norm{v_t}^2.
\label{eq:snapshot-descent}
\end{align}
The learning rate in \eqref{eq:snapshot-eta} satisfies
\[
L\eta=\frac1{2(1+\sqrt{(q+D)/b}+\tau)}\leq\frac12,
\]
and
\begin{align*}
2L^2\eta^2\left(\frac{q+D-1}{b}+\tau^2\right)
&\leq\frac{(q+D)/b+\tau^2}
{2(1+\sqrt{(q+D)/b}+\tau)^2}\\
&\leq\frac{(\sqrt{(q+D)/b}+\tau)^2}
{2(1+\sqrt{(q+D)/b}+\tau)^2}
\leq\frac12.
\end{align*}
Thus the bracket in \eqref{eq:snapshot-descent} is at least $1-1/2-1/2=0$. After dropping its nonpositive contribution and dividing by $T$, the uniform output rule gives
\[
\E\norm{\nabla f(x_R)}^2
=\frac1T\sum_{t=0}^{T-1}\E\norm{\nabla f(x_t)}^2
\leq\frac{2\Delta_f}{\eta T}.
\]
This proves \eqref{eq:snapshot-rate}.

There is one full gradient at initialization. The remaining snapshots are requested at positive multiples of $q$ smaller than $T$. Their total number, including initialization, is
$1+\lfloor(T-1)/q\rfloor=\lceil T/q\rceil$.
Every such request costs $n$ evaluations, whether or not its deadline occurs within the run. There are also $T$ two-point tasks with $b$ indices each, costing $2bT$ evaluations. This proves \eqref{eq:snapshot-cost}.
\end{proof}

\begin{proof}[Proof of Corollary~\ref{cor:snapshot}]
The choice $T\geq2\Delta_f/(\eta\epsilon^2)$ and \eqref{eq:snapshot-rate} give the target accuracy. With $q=n$ and $b=1$, we also have
\begin{align*}
T&\leq1+\frac{2\Delta_f}{\eta\epsilon^2}
=1+\frac{4L\Delta_f}{\epsilon^2}\bigl(1+\sqrt{n+D}+\tau\bigr)\\
&\leq1+\frac{4L\Delta_f}{\epsilon^2}
\bigl(1+\sqrt n+\sqrt D+\tau\bigr).
\end{align*}
Here $\sqrt{n+D}\leq\sqrt n+\sqrt D$ follows by squaring both nonnegative sides. Substitution into \eqref{eq:snapshot-cost} yields
\begin{align*}
n\left\lceil\frac Tn\right\rceil+2T
&\leq n+3T\\
&\leq n+3+\frac{12L\Delta_f}{\epsilon^2}
\bigl(1+\sqrt n+\sqrt D+\tau\bigr).
\end{align*}
Since $n\geq1$, we have $n+3\leq4n$ and $1\leq\sqrt n$. These two bounds give \eqref{eq:snapshot-complexity}.
\end{proof}

\section{The Effect of Outstanding Feedback}
\label{sec:outstanding}

The bounds based on maximum delay can be conservative when only a few corrections are missing. We first account for the number of outstanding corrections. For the stochastic method, we then use their total weight in the estimator, which also allows a sequence of permanently missing responses.

The maximum number of outstanding tasks and related concurrency quantities have been used in delayed learning and asynchronous SGD \citep{joulani2013delayed,koloskova2022sharper}. Our purpose here is to establish what these quantities imply for the two recursive variance-reduction estimators. The algorithms are unchanged.

\subsection{A bound for missing corrections}

The following inequality is deterministic. Therefore, it can be applied to every realized execution before any expectation is taken.

\begin{lemma}
\label{lem:weighted-missing}
For arbitrary vectors $\delta_j$, sets $\mathcal M_t\subseteq\{0,\ldots,t-1\}$, and nonnegative weights $w_{t,j}$,
\begin{align}
&\sum_{t=0}^{T-1}\norm{\sum_{j\in\mathcal M_t}w_{t,j}\delta_j}^2\nonumber\\
&\quad\leq\sum_{j=0}^{T-1}
\left[\sum_{t:j\in\mathcal M_t}w_{t,j}
\left(\sum_{i\in\mathcal M_t}w_{t,i}\right)\right]\norm{\delta_j}^2.
\label{eq:weighted-missing}
\end{align}
For any $p\geq0$, the following consequences hold:
\begin{enumerate}[label=(\roman*),topsep=3pt,itemsep=4pt]
\item If $|\mathcal M_t|\leq p$, $0\leq w_{t,j}\leq1$, and each task belongs to at most $\tau$ of the sets $\mathcal M_t$, the right-hand side of \eqref{eq:weighted-missing} is at most $p\tau\sum_j\norm{\delta_j}^2$.
\item For the stochastic weights $w_{t,j}=(1-\alpha)^{t-j-1}$, suppose that $\sum_{j\in\mathcal M_t}w_{t,j}\leq p$ at every round. Then the bound is $(p/\alpha)\sum_j\norm{\delta_j}^2$, without an upper bound on arrival times.
\item For the finite-sum method with $w_{t,j}=1$, epoch length $q$, and $|\mathcal M_t|\leq p$, the bound is $p(q-1)\sum_j\norm{\delta_j}^2$, without an upper bound on arrival times.
\end{enumerate}
\end{lemma}

\begin{proof}
The triangle inequality and scalar Cauchy--Schwarz give
\begin{align*}
\norm{\sum_{j\in\mathcal M_t}w_{t,j}\delta_j}^2
&\leq\left(\sum_{j\in\mathcal M_t}w_{t,j}\norm{\delta_j}\right)^2\\
&=\left(\sum_{j\in\mathcal M_t}
\sqrt{w_{t,j}}\,\sqrt{w_{t,j}}\norm{\delta_j}\right)^2\\
&\leq\left(\sum_{i\in\mathcal M_t}w_{t,i}\right)
\left(\sum_{j\in\mathcal M_t}w_{t,j}\norm{\delta_j}^2\right).
\end{align*}
Empty sums are zero. Summing over $t$ and changing the order of the two finite sums yields
\begin{align*}
&\sum_{t=0}^{T-1}\norm{\sum_{j\in\mathcal M_t}w_{t,j}\delta_j}^2\\
&\quad\leq\sum_{t=0}^{T-1}\sum_{j\in\mathcal M_t}
w_{t,j}\left(\sum_{i\in\mathcal M_t}w_{t,i}\right)\norm{\delta_j}^2\\
&\quad=\sum_{j=0}^{T-1}\left[
\sum_{\substack{0\leq t<T\\j\in\mathcal M_t}}
w_{t,j}\left(\sum_{i\in\mathcal M_t}w_{t,i}\right)
\right]\norm{\delta_j}^2.
\end{align*}
This is \eqref{eq:weighted-missing}; the coefficient for $j=T-1$ is zero.

For (i), the coefficient of each $\norm{\delta_j}^2$ satisfies
\begin{align*}
\sum_{t:j\in\mathcal M_t}w_{t,j}
\left(\sum_{i\in\mathcal M_t}w_{t,i}\right)
&\leq\sum_{t:j\in\mathcal M_t}1\cdot|\mathcal M_t|\\
&\leq p\sum_{t:j\in\mathcal M_t}1\leq p\tau.
\end{align*}
For (ii), use the assumed bound on the total missing weight and then sum a geometric series:
\begin{align*}
&\sum_{t:j\in\mathcal M_t}(1-\alpha)^{t-j-1}
\left(\sum_{i\in\mathcal M_t}(1-\alpha)^{t-i-1}\right)\\
&\quad\leq p\sum_{t=j+1}^{\infty}(1-\alpha)^{t-j-1}
=p\sum_{k=0}^{\infty}(1-\alpha)^k=\frac p\alpha.
\end{align*}
For (iii), let $r$ be the first round of the epoch containing task $j$. Then $r\leq j\leq r+q-1$, and
\[
j\in\mathcal M_t\quad\Longrightarrow\quad
j+1\leq t\leq r+q-1.
\]
There are at most $r+q-1-j\leq q-1$ such rounds. Since every weight is one,
\[
\sum_{t:j\in\mathcal M_t}\left(\sum_{i\in\mathcal M_t}1\right)
\leq p(r+q-1-j)\leq p(q-1).
\]
Substituting these three coefficient bounds into \eqref{eq:weighted-missing} proves the claims. Every inequality holds for the realized sets and vectors, even when the sets depend on the vectors.
\end{proof}

The bounds $p$ and $\tau$ used below are fixed upper bounds valid throughout the run, up to events of probability zero.

\Needspace{17\baselineskip}
\subsection{Few outstanding tasks under bounded delay}
\begin{corollary}
\label{cor:concurrency}
Suppose that \eqref{eq:delay} holds and $|\mathcal M_t|\leq p$ at every round. For the finite-sum method, $\mathcal M_t$ contains only current-epoch corrections. Theorem~\ref{thm:online} remains valid after replacing $\tau$ by $\sqrt{p\tau}$ in the choices of $\eta$, $\alpha$, and $T$, and then choosing $B$ by \eqref{eq:online-batch}. In Theorem~\ref{thm:finite}, replace the term $\min\{\tau,q-1\}$ in the denominator of the learning rate by
\[
\sqrt{\min\{p,\tau,q-1\}\min\{\tau,q-1\}}.
\]
With $q=n$ and $b=1$ for the finite-sum method, the two complexities become
\begin{align}
&O\!\left(1+\frac{\sigma^2}{\epsilon^2}
+\frac{L\Delta_f}{\epsilon^2}\left(1+\frac{\sigma}{\epsilon}+\sqrt{p\tau}\right)\right),
\label{eq:online-concurrency}\\
&O\!\left(n+\frac{L\Delta_f}{\epsilon^2}
\left(\sqrt n+\sqrt{\min\{p,\tau,n-1\}\min\{\tau,n-1\}}\right)\right),
\label{eq:finite-concurrency}
\end{align}
respectively.
\end{corollary}

\begin{proof}
For the stochastic method, a task $j$ can be missing only at rounds
$j+1,\ldots,j+\tau$. Lemma~\ref{lem:weighted-missing}(i) and \eqref{eq:online-packet} therefore give
\begin{align*}
&\sum_{t=0}^{T-1}\E\norm{\sum_{j\in\mathcal M_t}(1-\alpha)^{t-j-1}\delta_j}^2\\
&\quad\leq p\tau\sum_{j=0}^{T-1}\E\norm{\delta_j}^2\\
&\quad\leq2L^2p\tau\sum_{j=0}^{T-1}\E\norm{x_{j+1}-x_j}^2\\
&\qquad
+2\alpha^2p\tau\sum_{j=0}^{T-1}\E\norm{\nabla f(x_j)}^2
+2\alpha^2p\tau\sigma^2T.
\end{align*}
Combining this with \eqref{eq:online-decomposition}, \eqref{eq:virtual-error}, and \eqref{eq:descent-sum}, exactly as in \eqref{eq:online-master}, yields
\begin{align*}
&(1-4\alpha^2p\tau)\sum_{t=0}^{T-1}\E\norm{\nabla f(x_t)}^2\\
&\quad+\left[1-L\eta-4L^2\eta^2\left(\frac1\alpha+p\tau\right)\right]
\sum_{t=0}^{T-1}\E\norm{v_t}^2\\
&\qquad\leq\frac{2\Delta_f}{\eta}
+\frac{2\sigma^2}{B}\min\left\{T,\frac1\alpha\right\}
+4\sigma^2(\alpha+\alpha^2p\tau)T.
\end{align*}
This is \eqref{eq:online-master} with $\tau^2$ replaced by $p\tau$. In particular, using the parameters with $\sqrt{p\tau}$ in place of $\tau$ gives
\begin{align*}
4\alpha^2p\tau
&=\frac{p\tau}{64(1+\sigma/\epsilon)^2(1+\sigma/\epsilon+\sqrt{p\tau})^2}
\leq\frac1{64},\\
L\eta&\leq\frac1{16},\qquad
\frac{4L^2\eta^2}{\alpha}
=\frac{1+\sigma/\epsilon}{4(1+\sigma/\epsilon+\sqrt{p\tau})}\leq\frac14,\\
4L^2\eta^2p\tau
&=\frac{p\tau}{64(1+\sigma/\epsilon+\sqrt{p\tau})^2}\leq\frac1{64}.
\end{align*}
Thus the gradient coefficient is at least $63/64$, and the estimator coefficient is at least $43/64$. The same substitutions give
\begin{align*}
\frac{2\Delta_f}{\eta T}&\leq\frac{\epsilon^2}{2},
&\frac{2\sigma^2}{B}\min\left\{1,\frac1{\alpha T}\right\}
&\leq\frac{\epsilon^2}{8},\\
4\alpha\sigma^2&\leq\frac{\epsilon^2}{4},
&4\alpha^2p\tau\sigma^2&\leq\frac{\epsilon^2}{64}.
\end{align*}
Dropping the estimator term and dividing by $T$ therefore gives
\[
\frac1T\sum_{t=0}^{T-1}\E\norm{\nabla f(x_t)}^2
\leq\frac{\epsilon^2/2+\epsilon^2/8+\epsilon^2/4+\epsilon^2/64}{63/64}
=\frac{19}{21}\epsilon^2.
\]
Counting $B+2T$ evaluations proves \eqref{eq:online-concurrency}.

For the finite-sum method, at a round $t$ in an epoch beginning at $r$,
\[
\mathcal M_t\subseteq
\{\max\{r,t-\tau\},\ldots,t-1\},
\qquad |\mathcal M_t|\leq\min\{p,\tau,q-1\}.
\]
A task $j$ from this epoch can be missing only when
\[
j+1\leq t\leq\min\{j+\tau,r+q-1\}.
\]
The number of these rounds is at most $\min\{\tau,q-1\}$. Lemma~\ref{lem:weighted-missing} and \eqref{eq:finite-packet} now imply
\begin{align*}
&\sum_{t=0}^{T-1}\E\norm{\sum_{j\in\mathcal M_t}\delta_j}^2\\
&\quad\leq\min\{p,\tau,q-1\}\min\{\tau,q-1\}
\sum_{j=0}^{T-1}\E\norm{\delta_j}^2\\
&\quad\leq L^2\min\{p,\tau,q-1\}\min\{\tau,q-1\}
\sum_{j=0}^{T-1}\E\norm{x_{j+1}-x_j}^2.
\end{align*}
The complete-feedback bound \eqref{eq:finite-complete-sum} is unchanged. Consequently, the estimator coefficient after descent is
\[
1-L\eta-2L^2\eta^2
\left(\frac{q-1}{b}+\min\{p,\tau,q-1\}\min\{\tau,q-1\}\right).
\]
The stated learning rate makes this coefficient nonnegative: $L\eta\leq1/2$, and
\begin{align*}
&\frac{q-1}{b}+\min\{p,\tau,q-1\}\min\{\tau,q-1\}\\
&\quad\leq\left(1+\sqrt{q/b}
+\sqrt{\min\{p,\tau,q-1\}\min\{\tau,q-1\}}\right)^2
=\frac1{4L^2\eta^2}.
\end{align*}
It follows that $\E\norm{\nabla f(x_R)}^2\leq2\Delta_f/(\eta T)$. At $q=n,b=1$, the count is at most $n+3T$; substituting this learning rate and the prescribed $T$ proves \eqref{eq:finite-concurrency}.
\end{proof}

Under the delay assumption alone, $p=\tau$ recovers the original results. For example, suppose that one correction is delayed for $\tau\geq1$ rounds and all others return for the next update. Then $p=1$, and its contribution is controlled by $\sqrt\tau$. For the finite-sum method, the same comparison applies within an epoch long enough to contain the delayed interval. Related dependence on $\sqrt{p\tau}$ is known for asynchronous SGD \citep{koloskova2022sharper}.

\subsection{No maximum-delay assumption}

We now allow $A_j\in\{j+1,j+2,\ldots\}\cup\{\infty\}$, where $A_j=\infty$ means that the response never returns. The initialization batch and finite-sum full gradients must still be returned in full. Delay and non-return decisions may depend on already generated samples but not on future samples.

One possible condition is a bound on the number of missing responses:
\begin{equation}
|\mathcal M_t|\leq p\qquad\text{for every round }t.
\label{eq:pending-only}
\end{equation}
For the stochastic estimator, each old correction has a decreasing weight. We can therefore use the weaker condition
\begin{equation}
\sum_{j\in\mathcal M_t}(1-\alpha)^{t-j-1}\leq p
\qquad\text{for every round }t,
\label{eq:weighted-pending}
\end{equation}
where $\mathcal M_t$ includes every historical response that has not returned. Condition~\eqref{eq:pending-only} implies \eqref{eq:weighted-pending}, but \eqref{eq:weighted-pending} does not bound the number of permanently missing responses.

\begin{theorem}[Stochastic optimization with bounded missing weight]
\label{thm:online-pending}
Under Assumptions~\ref{ass:lower} and~\ref{ass:stochastic}, fix $p\geq0$ and $\epsilon>0$, and choose
\begin{align}
\eta&=\frac{1}{16L(1+p)(1+\sigma/\epsilon)},
&\alpha&=\frac{1}{16(1+p)(1+\sigma/\epsilon)^2},
\label{eq:pending-params-a}\\
T&=\max\left\{1,\left\lceil
\frac{128L\Delta_f(1+p)(1+\sigma/\epsilon)}{\epsilon^2}
\right\rceil\right\}.
\label{eq:pending-params-b}
\end{align}
Choose $B$ by \eqref{eq:online-batch} with these values of $\alpha$ and $T$. If \eqref{eq:weighted-pending} holds for this choice of $\alpha$, then Algorithm~\ref{alg:online} guarantees
\begin{equation}
\E\norm{\nabla f(x_R)}^2\leq\epsilon^2.
\label{eq:online-pending-target}
\end{equation}
The method issues at most $B+2T$ evaluations, giving the oracle complexity
\begin{equation}
O\!\left(1+\frac{\sigma^2}{\epsilon^2}
+\frac{L\Delta_f(1+p)}{\epsilon^2}\left(1+\frac{\sigma}{\epsilon}\right)\right).
\label{eq:online-pending-complexity}
\end{equation}
\end{theorem}

The order of these choices matters: $p$ is prescribed first, the algorithm then fixes $\alpha$, and the theorem applies to schedules satisfying \eqref{eq:weighted-pending} for that $\alpha$. A value estimated from the completed run is not used to set the learning rate.

\begin{proof}
Lemma~\ref{lem:weighted-missing}(ii) first removes the missing sets from the bound. Taking expectations and applying \eqref{eq:online-packet}, we obtain
\begin{align*}
&\sum_{t=0}^{T-1}\E\norm{\sum_{j\in\mathcal M_t}(1-\alpha)^{t-j-1}\delta_j}^2\\
&\quad\leq\frac{p}{\alpha}\sum_{j=0}^{T-1}\E\norm{\delta_j}^2\\
&\quad\leq\frac p\alpha\sum_{j=0}^{T-1}
\left(2L^2\E\norm{x_{j+1}-x_j}^2
+2\alpha^2\E\norm{\nabla f(x_j)}^2+2\alpha^2\sigma^2\right)\\
&\quad=\frac{2L^2p}{\alpha}\sum_{j=0}^{T-1}\E\norm{x_{j+1}-x_j}^2\\
&\qquad
+2\alpha p\sum_{j=0}^{T-1}\E\norm{\nabla f(x_j)}^2
+2\alpha p\sigma^2T.
\end{align*}
Using \eqref{eq:online-decomposition} and $\norm{a-b}^2\leq2\norm{a}^2+2\norm{b}^2$, and then \eqref{eq:virtual-error}, gives
\begin{align}
\sum_{t=0}^{T-1}\E\norm{v_t-\nabla f(x_t)}^2
\leq{}&2\sum_{t=0}^{T-1}\E\norm{u_t-\nabla f(x_t)}^2\nonumber\\
&+2\sum_{t=0}^{T-1}\E\norm{\sum_{j\in\mathcal M_t}(1-\alpha)^{t-j-1}\delta_j}^2\nonumber\\
\leq{}&\frac{2\sigma^2}{B}\min\left\{T,\frac1\alpha\right\}
+\frac{4L^2(1+p)}{\alpha}\sum_{t=0}^{T-1}\E\norm{x_{t+1}-x_t}^2\nonumber\\
&+4\alpha p\sum_{t=0}^{T-1}\E\norm{\nabla f(x_t)}^2
+4\alpha(1+p)\sigma^2T.
\label{eq:online-pending-error}
\end{align}
Substitute \eqref{eq:online-pending-error} into \eqref{eq:descent-sum}, use
$\norm{x_{t+1}-x_t}^2=\eta^2\norm{v_t}^2$, and move the gradient and estimator sums to the left. Then
\begin{align}
&(1-4\alpha p)\sum_{t=0}^{T-1}\E\norm{\nabla f(x_t)}^2\nonumber\\
&\quad+\left(1-L\eta-\frac{4L^2\eta^2(1+p)}{\alpha}\right)
\sum_{t=0}^{T-1}\E\norm{v_t}^2\nonumber\\
&\qquad\leq\frac{2\Delta_f}{\eta}
+\frac{2\sigma^2}{B}\min\left\{T,\frac1\alpha\right\}
+4\alpha(1+p)\sigma^2T.
\label{eq:online-pending-master}
\end{align}
The chosen parameters satisfy
\begin{align*}
4\alpha p
&=\frac{p}{4(1+p)(1+\sigma/\epsilon)^2}\leq\frac14,\\
L\eta&=\frac1{16(1+p)(1+\sigma/\epsilon)}\leq\frac1{16},\\
\frac{4L^2\eta^2(1+p)}{\alpha}
&=\frac{64(1+p)^2(1+\sigma/\epsilon)^2}
{256(1+p)^2(1+\sigma/\epsilon)^2}=\frac14.
\end{align*}
Thus the two coefficients in \eqref{eq:online-pending-master} are at least
\[
1-4\alpha p\geq\frac34,\qquad
1-L\eta-\frac{4L^2\eta^2(1+p)}{\alpha}
\geq1-\frac1{16}-\frac14=\frac{11}{16}.
\]
The choices of $T$ and $B$ imply
\begin{align*}
T&\geq\frac{128L\Delta_f(1+p)(1+\sigma/\epsilon)}{\epsilon^2}
=\frac{8\Delta_f}{\eta\epsilon^2},\\
B&\geq\frac{16\sigma^2}{\epsilon^2}
\min\left\{1,\frac1{\alpha T}\right\}.
\end{align*}
Consequently,
\begin{align*}
\frac{2\Delta_f}{\eta T}&\leq\frac{\epsilon^2}{4},\qquad
\frac{2\sigma^2}{B}\min\left\{1,\frac1{\alpha T}\right\}
\leq\frac{\epsilon^2}{8},\\
4\alpha(1+p)\sigma^2
&=\frac{\epsilon^2}{4}\left(\frac{\sigma/\epsilon}{1+\sigma/\epsilon}\right)^2
\leq\frac{\epsilon^2}{4}.
\end{align*}
The first two bounds also hold when $\Delta_f=0$ and $\sigma=0$, respectively, since the corresponding numerators are zero and $T,B\geq1$.
Dropping the nonnegative estimator term in \eqref{eq:online-pending-master} and dividing by $T$ gives
\[
\frac1T\sum_{t=0}^{T-1}\E\norm{\nabla f(x_t)}^2
\leq\frac{\epsilon^2/4+\epsilon^2/8+\epsilon^2/4}{3/4}
=\frac56\epsilon^2.
\]
Since $R$ is uniform and independent,
\[
\E\norm{\nabla f(x_R)}^2
=\sum_{t=0}^{T-1}\frac1T\E\norm{\nabla f(x_t)}^2
\leq\frac56\epsilon^2\leq\epsilon^2.
\]
Finally, $\max\{1,\lceil z\rceil\}\leq1+z$ for $z\geq0$, so
\begin{align*}
B&\leq1+\frac{16\sigma^2}{\epsilon^2},\\
T&\leq1+\frac{128L\Delta_f(1+p)(1+\sigma/\epsilon)}{\epsilon^2},\\
B+2T&\leq3+\frac{16\sigma^2}{\epsilon^2}
+\frac{256L\Delta_f(1+p)(1+\sigma/\epsilon)}{\epsilon^2}.
\end{align*}
This proves \eqref{eq:online-pending-complexity}, including all issued tasks.
\end{proof}

\begin{corollary}[A bound on the number of outstanding responses]
\label{cor:online-pending-count}
The conclusions of Theorem~\ref{thm:online-pending} hold if \eqref{eq:pending-only} holds for all historical missing corrections.
\end{corollary}
\begin{proof}
Since $j<t$ and $0<1-\alpha<1$,
\[
\sum_{j\in\mathcal M_t}(1-\alpha)^{t-j-1}
\leq\sum_{j\in\mathcal M_t}1
=|\mathcal M_t|\leq p.
\]
Thus \eqref{eq:weighted-pending} holds, and Theorem~\ref{thm:online-pending} applies with the same parameters.
\end{proof}

\paragraph{Permanently missing responses.}
As an example beyond the count condition, use $p=2$ in Theorem~\ref{thm:online-pending} and set
\[
m=\left\lceil\frac{\log 2}{\alpha}\right\rceil.
\]
Suppose that permanent losses have generation indices at least $m$ rounds apart, and every other response returns for the next update. At round $t$, write the missing indices in decreasing order as
$j_0>j_1>\cdots>j_{|\mathcal M_t|-1}$. Then
\[
j_0\leq t-1,\qquad j_k\leq j_0-km,
\qquad t-j_k-1\geq km.
\]
Since $0<1-\alpha<1$, increasing an exponent decreases its weight. Hence
\begin{align*}
\sum_{j\in\mathcal M_t}(1-\alpha)^{t-j-1}
&=\sum_{k=0}^{|\mathcal M_t|-1}(1-\alpha)^{t-j_k-1}\\
&\leq\sum_{k=0}^{\infty}(1-\alpha)^{km}
=\frac1{1-(1-\alpha)^m}\leq2,
\end{align*}
where the last step follows from
\[
\log(1-\alpha)\leq-\alpha,\qquad
(1-\alpha)^m\leq e^{-\alpha m}
\leq e^{-\log2}=\frac12.
\]
The empty missing set satisfies the bound as well. There is no upper bound on the total number of losses over an indefinitely long execution. Their spacing is accuracy dependent: at $p=2$,
\[
\alpha=\frac1{48(1+\sigma/\epsilon)^2},\qquad
\frac{\log2}{\alpha}\leq m<\frac{\log2}{\alpha}+1.
\]
Thus one loss every $m$ rounds has frequency $1/m=\Theta((1+\sigma/\epsilon)^{-2})$, which decreases as $\epsilon$ decreases when $\sigma>0$. Theorem~\ref{thm:online-pending} retains the stochastic variance-reduction order for these schedules. Every issued evaluation still counts; a missing response does not by itself mean that its physical computation has been cancelled.

\begin{theorem}[Finite sums with bounded outstanding corrections]
\label{thm:finite-pending}
Under Assumptions~\ref{ass:lower} and~\ref{ass:finite-smooth}, suppose that \eqref{eq:pending-only} holds for the current-epoch missing sets of Algorithm~\ref{alg:finite}. No maximum delay is assumed. Choose
\begin{equation}
\eta=\frac{1}{2L\bigl(1+\sqrt{q/b}+\sqrt{p(q-1)}\bigr)}.
\label{eq:finite-pending-eta}
\end{equation}
Then \eqref{eq:finite-rate} and \eqref{eq:finite-cost} hold. With $q=n$, $b=1$, and $T=\max\{1,\lceil2\Delta_f/(\eta\epsilon^2)\rceil\}$, the oracle complexity for \eqref{eq:target} is
\begin{equation}
O\!\left(n+\frac{L\Delta_f\sqrt{n(1+p)}}{\epsilon^2}\right).
\label{eq:finite-pending-complexity}
\end{equation}
\end{theorem}

\begin{proof}
By Lemma~\ref{lem:weighted-missing}(iii) and \eqref{eq:finite-packet},
\begin{align*}
\sum_{t=0}^{T-1}\E\norm{\sum_{j\in\mathcal M_t}\delta_j}^2
&\leq p(q-1)\sum_{j=0}^{T-1}\E\norm{\delta_j}^2\\
&\leq L^2p(q-1)\sum_{j=0}^{T-1}\E\norm{x_{j+1}-x_j}^2.
\end{align*}
The decomposition \eqref{eq:finite-decomposition} and the sampling-error estimate \eqref{eq:finite-complete-sum} do not use a delay bound. Applying $\norm{a-b}^2\leq2\norm{a}^2+2\norm{b}^2$ therefore gives
\begin{align}
\sum_{t=0}^{T-1}\E\norm{v_t-\nabla f(x_t)}^2
\leq{}&2\sum_{t=0}^{T-1}\E\norm{\sum_{j=r_t}^{t-1}\zeta_j}^2\nonumber\\
&+2\sum_{t=0}^{T-1}\E\norm{\sum_{j\in\mathcal M_t}\delta_j}^2\nonumber\\
\leq{}&2L^2\left(\frac{q-1}{b}+p(q-1)\right)
\sum_{t=0}^{T-1}\E\norm{x_{t+1}-x_t}^2.
\label{eq:finite-pending-error}
\end{align}
Substituting into \eqref{eq:descent-sum} and using $x_{t+1}-x_t=-\eta v_t$, we obtain
\begin{align*}
\sum_{t=0}^{T-1}\E\norm{\nabla f(x_t)}^2
\leq{}&\frac{2\Delta_f}{\eta}\\
&-\left[1-L\eta-2L^2\eta^2
\left(\frac{q-1}{b}+p(q-1)\right)\right]
\sum_{t=0}^{T-1}\E\norm{v_t}^2.
\end{align*}
For the learning rate \eqref{eq:finite-pending-eta},
\begin{align*}
L\eta&=\frac1{2(1+\sqrt{q/b}+\sqrt{p(q-1)})}\leq\frac12,\\
2L^2\eta^2\left(\frac{q-1}{b}+p(q-1)\right)
&\leq\frac{q/b+p(q-1)}{2(1+\sqrt{q/b}+\sqrt{p(q-1)})^2}
\leq\frac12.
\end{align*}
The last inequality uses $a^2+c^2\leq(1+a+c)^2$ for nonnegative $a,c$. The bracket is therefore nonnegative. Dropping its nonpositive contribution and dividing by $T$ proves
\[
\E\norm{\nabla f(x_R)}^2
=\frac1T\sum_{t=0}^{T-1}\E\norm{\nabla f(x_t)}^2
\leq\frac{2\Delta_f}{\eta T},
\]
which is \eqref{eq:finite-rate}. Refreshes still occur at rounds $0,q,2q,\ldots<T$, and at most $T$ correction tasks are issued. Their total cost is at most $n\lceil T/q\rceil+2bT$, proving \eqref{eq:finite-cost}.

For $q=n,b=1$, each of $1$, $\sqrt n$, and $\sqrt{p(n-1)}$ is at most $\sqrt{n(1+p)}$. It follows that
\begin{align*}
\eta^{-1}&=2L(1+\sqrt n+\sqrt{p(n-1)})
\leq6L\sqrt{n(1+p)},\\
T&\leq1+\frac{2\Delta_f}{\eta\epsilon^2}
\leq1+\frac{12L\Delta_f\sqrt{n(1+p)}}{\epsilon^2},\\
n\lceil T/n\rceil+2T&\leq n+3T
\leq n+3+\frac{36L\Delta_f\sqrt{n(1+p)}}{\epsilon^2}.
\end{align*}
The prescribed $T$ gives $2\Delta_f/(\eta T)\leq\epsilon^2$, and $n\geq1$ absorbs the constant three into the stated complexity.
\end{proof}

\paragraph{Interpretation.}
For a fixed $p$, Theorems~\ref{thm:online-pending} and~\ref{thm:finite-pending} retain their respective variance-reduction orders without a maximum age for a correction. The stochastic result bounds the remaining weight of missing responses, whereas the finite-sum result bounds their count within the current epoch. Exponential weights gradually remove the influence of an old stochastic correction; an exact refresh removes the influence of an old finite-sum correction. Both arguments allow sample-dependent missingness. Appendix~\ref{app:presampled} gives a sharper dependence on the outstanding count when each arrival time is fixed before drawing its task's sample.

\bibliographystyle{plainnat}
\begingroup
\footnotesize
\setlength{\bibsep}{4pt}
\bibliography{references}
\endgroup
\clearpage
\appendix
\section{Arrival Times Fixed Before Sampling}
\label{app:presampled}

The analysis in the main text allows a task's completion time to depend on its sampled gradient. In this appendix, we impose a stronger timing condition and improve the dependence on the number of outstanding corrections. Algorithm~\ref{alg:online} is unchanged.

\begin{assumption}[Arrival time chosen before the sample]
\label{ass:presampled}
For each task $j$, after its point pair has been chosen, the arrival time $A_j$ is fixed before drawing $\xi_j$. The arrival time may depend on earlier samples and independent scheduling randomness, and may equal infinity. It is not revised later. The sample $\xi_j$ has the prescribed oracle distribution after the point pair and $A_j$ have been fixed.
\end{assumption}

This condition includes a schedule fixed in advance. It also permits later arrival times to depend on earlier samples. No maximum delay is required.

\begin{theorem}
\label{thm:presampled}
Under Assumptions~\ref{ass:lower}, \ref{ass:stochastic}, and~\ref{ass:presampled}, suppose that \eqref{eq:pending-only} holds for all historical missing corrections, with a fixed $p\geq0$. Given $\epsilon>0$, choose
\begin{align}
\alpha&=\frac{1}{16\max\{(1+\sigma/\epsilon)^2,1+p\}},
&\eta&=\frac{\sqrt\alpha}{4L\sqrt{1+p}},
\label{eq:presampled-params}\\
T&=\max\left\{1,\left\lceil\frac{8\Delta_f}{\eta\epsilon^2}\right\rceil\right\},
\label{eq:presampled-horizon}
\end{align}
and choose $B$ by \eqref{eq:online-batch}. Then Algorithm~\ref{alg:online} satisfies
\[
\E\norm{\nabla f(x_R)}^2\leq\epsilon^2.
\]
It issues at most $B+2T$ evaluations, with complexity
\begin{equation}
O\!\left(
1+\frac{\sigma^2}{\epsilon^2}
+\frac{L\Delta_f}{\epsilon^2}
\max\left\{\left(1+\frac{\sigma}{\epsilon}\right)\sqrt{1+p},1+p\right\}
\right).
\label{eq:presampled-complexity}
\end{equation}
\end{theorem}

For comparison, Corollary~\ref{cor:online-pending-count} contains the factor $(1+p)(1+\sigma/\epsilon)$. Under Assumption~\ref{ass:presampled}, the stochastic term grows as $\sqrt{1+p}$, while the deterministic term grows as $1+p$.

\subsection{Estimation error}

Recall the centered innovation $\zeta_j$ in \eqref{eq:online-noise-definition}. Separate the received innovations from the missing population corrections by defining
\begin{equation}
a_t=(1-\alpha)^t\bigl(v_0-\nabla f(x_0)\bigr)
+\sum_{j=0}^{t-1}\ind_{\{A_j\leq t\}}(1-\alpha)^{t-j-1}\zeta_j.
\label{eq:presampled-noise}
\end{equation}
The population corrections telescope:
\begin{align*}
&\sum_{j=0}^{t-1}(1-\alpha)^{t-j-1}
\bigl(\nabla f(x_{j+1})-(1-\alpha)\nabla f(x_j)\bigr)\\
&\quad=\sum_{j=1}^{t}(1-\alpha)^{t-j}\nabla f(x_j)
-\sum_{j=0}^{t-1}(1-\alpha)^{t-j}\nabla f(x_j)\\
&\quad=\nabla f(x_t)-(1-\alpha)^t\nabla f(x_0).
\end{align*}
Subtracting this identity from \eqref{eq:online-estimator} gives
\begin{align}
v_t-\nabla f(x_t)
={}&a_t-\sum_{j\in\mathcal M_t}(1-\alpha)^{t-j-1}
\bigl(\nabla f(x_{j+1})-(1-\alpha)\nabla f(x_j)\bigr).
\label{eq:presampled-decomposition}
\end{align}

\begin{lemma}
\label{lem:presampled-error}
Under the assumptions of Theorem~\ref{thm:presampled}, for any $0<\alpha<1$,
\begin{align}
\sum_{t=0}^{T-1}\E\norm{v_t-\nabla f(x_t)}^2
\leq{}&\frac{2\sigma^2}{B}\min\left\{T,\frac1\alpha\right\}
+\frac{4L^2(1+p)}{\alpha}
\sum_{t=0}^{T-1}\E\norm{x_{t+1}-x_t}^2\nonumber\\
&+4\alpha p\sum_{t=0}^{T-1}\E\norm{\nabla f(x_t)}^2
+4\alpha\sigma^2T.
\label{eq:presampled-error}
\end{align}
\end{lemma}

\begin{proof}
For each $j<t$, the initialization error, $A_j$, and the point pair of task $j$ are fixed before drawing $\xi_j$. For $i<j$, the same holds for $A_i$ and $\zeta_i$. Averaging the fresh sample $\xi_j$ first gives
\begin{align*}
\E\left[\ind_{\{A_i\leq t\}}\ind_{\{A_j\leq t\}}
\inner{\zeta_i}{\zeta_j}\right]&=0\quad(i<j),\\
\E\left[\ind_{\{A_j\leq t\}}
\inner{v_0-\nabla f(x_0)}{\zeta_j}\right]&=0.
\end{align*}
Thus, expanding \eqref{eq:presampled-noise}, we obtain
\begin{align*}
\E\norm{a_t}^2
={}&(1-\alpha)^{2t}\E\norm{v_0-\nabla f(x_0)}^2\\
&+\sum_{j=0}^{t-1}(1-\alpha)^{2(t-j-1)}
\E\left[\ind_{\{A_j\leq t\}}\norm{\zeta_j}^2\right]\\
&+2\sum_{j=0}^{t-1}(1-\alpha)^{2t-j-1}
\underbrace{\E\left[\ind_{\{A_j\leq t\}}
\inner{v_0-\nabla f(x_0)}{\zeta_j}\right]}_{=0}\\
&+2\sum_{0\leq i<j<t}(1-\alpha)^{2t-i-j-2}
\underbrace{\E\left[\ind_{\{A_i\leq t\}}\ind_{\{A_j\leq t\}}
\inner{\zeta_i}{\zeta_j}\right]}_{=0}.
\end{align*}
Using \eqref{eq:online-noise} and $\ind_{\{A_j\leq t\}}\leq1$ gives
\begin{align*}
\E\norm{a_t}^2
\leq{}&(1-\alpha)^{2t}\frac{\sigma^2}{B}\\
&+\sum_{j=0}^{t-1}(1-\alpha)^{2(t-j-1)}
\left(2L^2\E\norm{x_{j+1}-x_j}^2+2\alpha^2\sigma^2\right),
\end{align*}
where $\E\norm{v_0-\nabla f(x_0)}^2\leq\sigma^2/B$ follows from the initialization batch. The two geometric sums satisfy
\[
\sum_{t=0}^{T-1}(1-\alpha)^{2t}\leq\min\{T,1/\alpha\},
\qquad
\sum_{t=j+1}^{T-1}(1-\alpha)^{2(t-j-1)}\leq\frac1\alpha.
\]
Changing the order of summation now yields
\begin{align*}
\sum_{t=0}^{T-1}\E\norm{a_t}^2
\leq{}&\frac{\sigma^2}{B}\min\{T,1/\alpha\}\\
&+\sum_{j=0}^{T-2}
\left(2L^2\E\norm{x_{j+1}-x_j}^2+2\alpha^2\sigma^2\right)
\sum_{t=j+1}^{T-1}(1-\alpha)^{2(t-j-1)}\\
\leq{}&\frac{\sigma^2}{B}\min\{T,1/\alpha\}
+\frac{2L^2}{\alpha}\sum_{j=0}^{T-2}\E\norm{x_{j+1}-x_j}^2
+2\alpha\sigma^2(T-1).
\end{align*}
Extending the nonnegative sums to $T$ terms gives
\begin{equation}
\sum_{t=0}^{T-1}\E\norm{a_t}^2
\leq\frac{\sigma^2}{B}\min\left\{T,\frac1\alpha\right\}
+\frac{2L^2}{\alpha}\sum_{t=0}^{T-1}\E\norm{x_{t+1}-x_t}^2
+2\alpha\sigma^2T.
\label{eq:presampled-noise-sum}
\end{equation}

Next consider the missing population corrections. Their squared norms satisfy
\begin{align*}
\norm{\nabla f(x_{j+1})-(1-\alpha)\nabla f(x_j)}^2
&=\norm{\nabla f(x_{j+1})-\nabla f(x_j)+\alpha\nabla f(x_j)}^2\\
&\leq2L^2\norm{x_{j+1}-x_j}^2+2\alpha^2\norm{\nabla f(x_j)}^2.
\end{align*}
Since $\sum_{j\in\mathcal M_t}(1-\alpha)^{t-j-1}\leq|\mathcal M_t|\leq p$, weighted Cauchy--Schwarz gives
\begin{align}
&\sum_{t=0}^{T-1}\E\norm{
\sum_{j\in\mathcal M_t}(1-\alpha)^{t-j-1}
\bigl(\nabla f(x_{j+1})-(1-\alpha)\nabla f(x_j)\bigr)
}^2\nonumber\\
&\quad\leq p\sum_{t=0}^{T-1}\sum_{j=0}^{t-1}(1-\alpha)^{t-j-1}
\E\!\left[\ind_{\{j\in\mathcal M_t\}}
\norm{\nabla f(x_{j+1})-(1-\alpha)\nabla f(x_j)}^2\right]\nonumber\\
&\quad\leq p\sum_{j=0}^{T-2}
\E\norm{\nabla f(x_{j+1})-(1-\alpha)\nabla f(x_j)}^2
\sum_{t=j+1}^{T-1}(1-\alpha)^{t-j-1}\nonumber\\
&\quad\leq\frac{p}{\alpha}\sum_{j=0}^{T-1}
\E\norm{\nabla f(x_{j+1})-(1-\alpha)\nabla f(x_j)}^2\nonumber\\
&\quad\leq\frac{2L^2p}{\alpha}\sum_{j=0}^{T-1}\E\norm{x_{j+1}-x_j}^2
+2\alpha p\sum_{j=0}^{T-1}\E\norm{\nabla f(x_j)}^2.
\label{eq:presampled-population-sum}
\end{align}
The second inequality drops the indicator and exchanges the two sums; the third uses $\sum_{t=j+1}^{\infty}(1-\alpha)^{t-j-1}=1/\alpha$.
Finally, \eqref{eq:presampled-decomposition} and $\norm{a-b}^2\leq2\norm a^2+2\norm b^2$ imply
\begin{align*}
\sum_{t=0}^{T-1}\E\norm{v_t-\nabla f(x_t)}^2
\leq{}&\frac{2\sigma^2}{B}\min\{T,1/\alpha\}
+\left(\frac{4L^2}{\alpha}+\frac{4L^2p}{\alpha}\right)
\sum_{t=0}^{T-1}\E\norm{x_{t+1}-x_t}^2\nonumber\\
&+4\alpha p\sum_{t=0}^{T-1}\E\norm{\nabla f(x_t)}^2
+4\alpha\sigma^2T,
\end{align*}
which is \eqref{eq:presampled-error}.
\end{proof}

\subsection{Proof of Theorem~\ref{thm:presampled}}

Substitute \eqref{eq:presampled-error} and
$\norm{x_{t+1}-x_t}^2=\eta^2\norm{v_t}^2$ into \eqref{eq:descent-sum}. Moving the gradient and estimator sums to the left gives
\begin{align}
&(1-4\alpha p)\sum_{t=0}^{T-1}\E\norm{\nabla f(x_t)}^2\nonumber\\
&\quad+\left(1-L\eta-\frac{4L^2\eta^2(1+p)}{\alpha}\right)
\sum_{t=0}^{T-1}\E\norm{v_t}^2\nonumber\\
&\qquad\leq\frac{2\Delta_f}{\eta}
+\frac{2\sigma^2}{B}\min\left\{T,\frac1\alpha\right\}
+4\alpha\sigma^2T.
\label{eq:presampled-master}
\end{align}
By \eqref{eq:presampled-params}, $\alpha\leq1/16$ and
\begin{align*}
4\alpha p
&=\frac{p}{4\max\{(1+\sigma/\epsilon)^2,1+p\}}
\leq\frac{p}{4(1+p)}\leq\frac14,\\
L\eta&=\frac{\sqrt\alpha}{4\sqrt{1+p}}\leq\frac1{16},\\
\frac{4L^2\eta^2(1+p)}{\alpha}
&=\frac{4L^2(1+p)}{\alpha}\cdot\frac{\alpha}{16L^2(1+p)}=\frac14.
\end{align*}
Hence the two coefficients on the left of \eqref{eq:presampled-master} satisfy
\[
1-4\alpha p\geq\frac34,
\qquad
1-L\eta-\frac{4L^2\eta^2(1+p)}{\alpha}
\geq1-\frac1{16}-\frac14=\frac{11}{16}>0.
\]
The definitions of $T$ and $B$ give
\[
\frac1T\leq\frac{\eta\epsilon^2}{8\Delta_f}\quad(\Delta_f>0),
\qquad
\frac{\sigma^2}{B}\min\left\{1,\frac1{\alpha T}\right\}\leq\frac{\epsilon^2}{16}.
\]
Therefore,
\begin{align*}
\frac{2\Delta_f}{\eta T}&\leq\frac{\epsilon^2}{4},\\
\frac{2\sigma^2}{B}\min\left\{1,\frac1{\alpha T}\right\}
&\leq\frac{\epsilon^2}{8},\\
4\alpha\sigma^2
&\leq\frac{\sigma^2}{4(1+\sigma/\epsilon)^2}
\leq\frac{\epsilon^2}{4}.
\end{align*}
The first bound also holds when $\Delta_f=0$, since its left-hand side is zero.
Dropping the nonnegative estimator term in \eqref{eq:presampled-master} gives
\[
\frac1T\sum_{t=0}^{T-1}\E\norm{\nabla f(x_t)}^2
\leq\frac{\epsilon^2/4+\epsilon^2/8+\epsilon^2/4}{3/4}
=\frac56\epsilon^2.
\]
For the uniform output,
\[
\E\norm{\nabla f(x_R)}^2
=\sum_{t=0}^{T-1}\frac1T\E\norm{\nabla f(x_t)}^2
\leq\frac56\epsilon^2\leq\epsilon^2.
\]
For the evaluation count,
\begin{align*}
\frac1\eta
&=16L\sqrt{(1+p)\max\{(1+\sigma/\epsilon)^2,1+p\}}\\
&=16L\max\left\{
\left(1+\frac\sigma\epsilon\right)\sqrt{1+p},1+p
\right\}.
\end{align*}
Since $B\leq1+16\sigma^2/\epsilon^2$ and $T\leq1+8\Delta_f/(\eta\epsilon^2)$,
\begin{align*}
B+2T
&\leq3+\frac{16\sigma^2}{\epsilon^2}+\frac{16\Delta_f}{\eta\epsilon^2}\\
&=3+\frac{16\sigma^2}{\epsilon^2}
+\frac{256L\Delta_f}{\epsilon^2}
\max\left\{\left(1+\frac\sigma\epsilon\right)\sqrt{1+p},1+p\right\},
\end{align*}
which proves \eqref{eq:presampled-complexity}. \qed

\paragraph{Why the timing condition is used.}
Let the stochastic gradient have additive noise $Z_j\in\{-\sigma,\sigma\}$, with equal probabilities. Then $\zeta_j=\alpha Z_j$. If $A_j=j+1$ for $Z_j=-\sigma$ and $A_j=j+2$ for $Z_j=\sigma$, we have
\[
\E\zeta_j=0,
\qquad
\E\left[\ind_{\{A_j\leq j+1\}}\zeta_j\right]
=\tfrac12(-\alpha\sigma)+\tfrac12\cdot0
=-\frac{\alpha\sigma}{2}.
\]
Thus an arrival indicator depending on its own sample need not preserve zero mean. The main-text analysis instead removes these indicators through a deterministic bound on the missing corrections.

\end{document}